\documentclass[11pt,a4paper]{article}
\usepackage[margin=27mm]{geometry}
\usepackage{amsmath,amssymb,amsthm,mathtools}
\usepackage[T1]{fontenc}
\usepackage{tikz}
\usetikzlibrary{arrows.meta,calc,positioning,patterns,decorations.pathreplacing,backgrounds}
\usepackage[font=small,labelfont=bf]{caption}
\usepackage[hidelinks]{hyperref}
\usepackage{enumitem}
\definecolor{coreblue}{HTML}{DDEAF2}
\definecolor{tailorange}{HTML}{E4A35A}
\definecolor{tailgreen}{HTML}{76A99B}
\definecolor{ink}{HTML}{26333D}
\definecolor{wallblue}{HTML}{23689A}
\definecolor{pathred}{HTML}{AD3D49}
\definecolor{safeolive}{HTML}{297365}
\definecolor{ghost}{HTML}{E8E8E8}
\tikzset{cell/.style={draw=ink!65,line width=.35pt},
 marker/.style={draw=ink,fill=white,inner sep=1.5pt,font=\small},
 wall/.style={draw=wallblue,line width=1.3pt},
 slack/.style={draw=pathred,line width=1.1pt,densely dashed},
 axis/.style={draw=ink!70,-{Stealth[length=2mm]}},
 paneltitle/.style={font=\bfseries\large,anchor=west},
 every node/.style={text=ink}}
\newtheorem{theorem}{Theorem}[section]
\newtheorem{lemma}[theorem]{Lemma}
\newtheorem{proposition}[theorem]{Proposition}
\newtheorem{corollary}[theorem]{Corollary}
\theoremstyle{definition}
\newtheorem{definition}[theorem]{Definition}
\theoremstyle{remark}
\newcommand{\NE}{\mathrm{NE}}
\newcommand{\NW}{\mathrm{NW}}

\title{Enumeration of $k$-convex polyominoes\\by semi-perimeter}
\author{Thomas Prellberg\thanks{Corresponding author.}\\[0.5ex]
\small School of Mathematical Sciences\\
\small Queen Mary University of London\\
\small Mile End Road, London E1 4NS, United Kingdom\\
\small\href{mailto:t.prellberg@qmul.ac.uk}{\texttt{t.prellberg@qmul.ac.uk}}}
\date{}
\hypersetup{pdftitle={Enumeration of k-convex polyominoes by semi-perimeter},
 pdfauthor={Thomas Prellberg}}
\begin{document}\maketitle
\begin{abstract}
We prove the Conway--Guttmann formula for the generating function of
$k$-convex polyominoes by semi-perimeter, for every $k\ge2$.  A decomposition
according to the two diagonal directions reduces the enumeration to convex
polyominoes satisfying a one-bend condition in one diagonal direction.
Splitting these polyominoes at the right endpoint of the top row leaves an
initial part followed by columns whose vertical intervals decrease by
inclusion.  Ordered forests describe the initial part and retain the row
information needed to determine when column extensions first exceed the
prescribed bound.  A first-passage calculation, together with a sum over
initial parts already above the bound, gives the required generating
function.  The argument includes the case $k=2$.
\end{abstract}


\section{Introduction}
We identify the unit square $[a,a+1]\times[b,b+1]$ with its integer
coordinates $(a,b)$.  A \emph{polyomino} is a finite nonempty set of such
cells that is connected by shared edges, considered up to translation.
Its width and height are the numbers of columns and rows that it meets.
Its \emph{perimeter} is the number of exposed unit edges, and half this
number is its \emph{semi-perimeter}, denoted $\operatorname{sp}(P)$.
It is \emph{convex} when every nonempty row and column is an interval of
cells.  For a convex polyomino the semi-perimeter equals width plus height:
each row contributes two vertical boundary edges and each column two
horizontal boundary edges.  Its \emph{bounding rectangle} is the smallest
rectangle with sides parallel to the coordinate axes that contains it.  Unless stated otherwise, $[a,b]$ denotes the integer interval
$\{a,a+1,\ldots,b\}$ when it describes rows or columns.

A path in $P$ is a finite sequence of its cells in which successive
cells share an edge.  It is \emph{monotone} when its steps use at most two coordinate
directions, which are not opposite.  A \emph{bend} is a change of step
direction; a path consisting of a single cell has no bends.  We write
$E,N,S,W$ for unit steps in the four directions.  Juxtaposition denotes
concatenation, so $N^aE^b$ means $a$ north steps followed by $b$ east steps;
an exponent zero contributes no steps.

Every convex polyomino contains a monotone path between any two cells.
Indeed, a shortest path cannot leave a row and subsequently return to it,
since row convexity would give a shorter horizontal replacement.  The same
argument applies to columns.

A convex polyomino is \emph{$k$-convex} if every pair of its cells can
be joined inside it by a monotone path with at most $k$ bends.  Bounding
this number refines convexity by a geometric measure of the paths needed
to connect cells.

Delest and Viennot~\cite{DV} obtained the generating function of the
unrestricted convex class.  The bound on the number of bends leads to the
hierarchy introduced by Castiglione and Restivo~\cite{CR}; its formulation
and early enumerative questions are discussed in~\cite{CFRR}.  The first
two nontrivial bounds already give different kinds of generating functions.
The $1$-convex polyominoes, also called $L$-convex polyominoes, have a rational
semi-perimeter generating function, obtained by Castiglione, Frosini,
Restivo and Rinaldi~\cite{CFRR}.  The $2$-convex or $Z$-convex class was
enumerated by Duchi, Rinaldi and Schaeffer~\cite{DRS}; its generating
function is algebraic.
Micheli and Rossin~\cite{MR} studied asymptotic lower bounds through
geometrically restricted subclasses of prescribed convexity degree.

A related family admits a formula for every degree bound.  A convex
polyomino is directed if all its cells can be reached from its southwest
corner cell by north and east steps.  Boussicault, Rinaldi and
Socci~\cite{BRS} represented these polyominoes by two ordered forests and
a lattice path, and obtained rational generating functions for the directed
$k$-convex class in terms of Fibonacci polynomials.  For the full convex
class, Brocchi, Castiglione and Massazza~\cite{BCM} studied how the two
diagonal directions determine the degree of convexity.  Beaton and
Rinaldi~\cite{BR} subsequently used these directional degrees to study
ascending polyominoes and to refine the enumeration of the $Z$-convex class.
The directional degrees and the term ascending are defined below.

Conway and Guttmann~\cite{CG} conjectured an explicit semi-perimeter
generating function for the full $k$-convex class.  We prove their formula
for every integer $k\ge2$.  The argument combines the directional
geometric decomposition with the forest representation, retaining the order
of the rows that meet the leftmost column.  This additional information
determines when successively appended columns first increase the degree
beyond the prescribed bound.  The resulting calculation treats the
contribution of the initial part and the first degree increase separately,
then sums both contributions explicitly.  The outline below explains the
geometric constructions before their detailed development.

\subsection{Statement of the result}

An ordered pair $a=(a_x,a_y),b=(b_x,b_y)$ is
\emph{northeast-comparable} when $a_x\le b_x$ and $a_y\le b_y$,
and \emph{northwest-comparable} when $a_x\ge b_x$ and $a_y\le b_y$.
For such a pair, $d_{\NE}(a,b)$ or $d_{\NW}(a,b)$ denotes the least
number of bends of a path inside $P$ in the indicated directions
from $a$ to $b$.
Their maxima over the corresponding ordered pairs in $P$ are the
\emph{directional degrees} $D_{\NE}(P)$ and $D_{\NW}(P)$.
The pair $a=b$ is allowed and has distance zero.  Path reversal shows that
\[
 \deg P=\max\{D_{\NE}(P),D_{\NW}(P)\}
\]
is the greatest minimum number of bends needed to join any two cells
by a monotone path.  Thus $P$ is $k$-convex precisely when $\deg P\le k$.
We call $P$ \emph{ascending} when $D_{\NW}(P)\le1$.

Write $Q_k(x)$ for the generating function of $k$-convex polyominoes by
semi-perimeter.  All series and substitutions below are formal; a square root
with constant term one always means that branch.  Define polynomials
\[
 d_0=1,\qquad d_1=1-2x,\qquad
 d_j=(1-2x)d_{j-1}-x^2d_{j-2}\quad(j\ge2),
\]
\[
 r_0=2,\qquad r_1=1,\qquad r_j=r_{j-1}-xr_{j-2}\quad(j\ge2),
\]
and set
\[
 p_k=(1-6x+7x^2+4x^3)d_k+x^2(1-4x+8x^2)d_{k-1}-2x^{k+2}.
\]
\begin{theorem}\label{thm:main}
For every integer $k\ge2$,
\begin{equation}\label{eq:main}
 \boxed{Q_k(x)=\frac{x^2p_k}{(1-4x)^2d_k}
                   -\frac{x^4r_k^2}{(1-4x)^{3/2}d_k}.}
\end{equation}
\end{theorem}
Theorem~\ref{thm:main} proves the polynomial form of the conjecture of
Conway and Guttmann~\cite[Section~4.1, pp.~79--80]{CG}.  It gives one
expression valid for every $k\ge2$; the case $k=2$ is included in the
same geometric decomposition.

\subsection{Outline of the proof}

Write $C(x)$ for the series of all convex polyominoes, $A(x)$ for
ascending polyominoes, and $A_k(x)$ for ascending polyominoes with
$D_{\NE}\le k$.  Each polyomino contributes $x^{\operatorname{sp}(P)}$.
We shall prove
\[
 Q_k=C-2(A-A_k)
\]
and obtain an explicit expression for $A-A_k$.  The reduction, proved in
Section~\ref{sec:directions}, uses the fact that a convex polyomino outside
the $k$-convex class has exactly one directional degree greater than $k$;
the other is at most one.  Reflection exchanges the two possibilities.
It therefore remains to count the ascending polyominoes whose northeast
degree exceeds $k$.

The geometric decomposition begins at the right endpoint of the top row.
The columns up to this point form the \emph{core}; the columns to its right
form the \emph{tail}, and their vertical intervals decrease by inclusion.
The core contains the northeast corner of its own bounding rectangle,
although it need not be a rectangle.  When the target cell of a northeast
path lies in the tail, only columns no farther right than the target and
rows no higher than the target can be used.  This leads to the \emph{shadows}: prefixes of
columns truncated above the top of their last column.  Section~\ref{pc:core-section}
shows that the northeast degree is the maximum of the degrees of these
shadows, including the untruncated core as the initial shadow, and is at
most one greater than the degree of the core.

After a half-turn, each core is a directed convex polyomino: every cell
can be reached from its southwest corner by north and east steps.
Section~\ref{pf:forest-section} uses the representation of Boussicault,
Rinaldi and Socci.  The directed polyomino is completed to one bounded by
two north--east paths; two ordered forests describe this completion, and a
separating path specifies the cells removed from its northeast corner.
The forest heights give the degree of the completion, and the position
of the separating path decides whether removing those cells lowers the
degree by one.  Keeping track also of the rows meeting the leftmost column
gives the information needed when that column is shortened or replaced.

The two ends of each nested tail interval move inwards.  Their cumulative
movements give two nonnegative coordinates.  Sections~\ref{sec:truncation}
and~\ref{sec:fixed-boundary} express the degree of each shadow in these
coordinates and identify a fixed monotone boundary, determined by the
original core, between extensions that satisfy the bound and those that
exceed it.  Thus the first increase beyond $k$ becomes the first passage
of a lattice path to the forbidden side of this boundary.  Before the
endpoint, a contact with the boundary is safe when the next boundary step
is east, and forbidden when it is north.

There are consequently two contributions to $A-A_k$.  Either the core
already has degree greater than $k$, or it has degree $k$ and the tail
has a first forbidden shadow.  Section~\ref{sec:separation} makes this
partition precise.  In Section~\ref{pe:enumeration}, the first contribution
is summed by a refinement of the forest series and a telescoping identity;
the second is obtained by counting the paths up to their first passage
and then adjoining the remaining nested columns.  The degree-two argument
in Section~\ref{sec:truncation} ensures that this decomposition also applies
at $k=2$.  Finally, Section~\ref{sec:conclusion} combines the two contributions
with the convex-polyomino series of Delest and Viennot~\cite{DV} to obtain
Theorem~\ref{thm:main}.

\subsection*{Research process}
ChatGPT (OpenAI) was used for mathematical exploration, and Claude
(Anthropic) was used to obtain critical feedback. These interactive tools
supported the research process. The results are established by the
mathematical constructions and proofs given in this article.
The explanatory figures were prepared in TikZ with assistance from
ChatGPT (OpenAI).

\section{Convexity and the two directional degrees}\label{sec:directions}

Write row $i$, with heights increasing from bottom to top, as $[L_i,R_i]$.
Two row intervals are \emph{nested} if either contains the other; equality
is allowed.  More generally, a sequence of intervals is nested when each
interval contains the next.

\begin{lemma}\label{pc:endpoint-shapes}
For a convex polyomino, $L_i$ first weakly decreases and then weakly increases,
whereas $R_i$ first weakly increases and then weakly decreases.
\end{lemma}
\begin{proof}
There is no strict peak $L_j>\max(L_i,L_\ell)$ with $i<j<\ell$.
Indeed put $x=L_j-1$ and take the nearest rows below and above $j$ whose left
endpoints are at most $x$.  Each overlaps its adjacent row towards $j$, whose
left endpoint exceeds $x$, and hence each contains column $x$.  Row $j$
does not, contradicting column convexity.  The absence of such peaks gives
the asserted valley shape.  Horizontal reflection proves the assertion for
$R_i$.
\end{proof}

Call a pair $i<j$ northeast shifted if $L_i<L_j$ and $R_i<R_j$, and northwest
shifted if $L_i>L_j$ and $R_i>R_j$.

The row criterion below is Beaton and Rinaldi's
Proposition~2.4~\cite{BR}.  The restriction on the two directional degrees
is due to Brocchi, Castiglione and Massazza~\cite[Proposition~15 and the
following discussion]{BCM}.  We give a proof of the precise form used here.
\begin{lemma}\label{pc:directional-structure}
A convex polyomino has $D_{\mathrm{NW}}\le1$ precisely when every pair of row intervals
is nested or northeast shifted.  Moreover, if both directional degrees exceed
one, both are at most two.
\end{lemma}
\begin{proof}
For northwest-comparable $a=(u,i)$ and $b=(v,j)$, where $i<j$ and $u\ge v$,
a path with at most one bend exists exactly when one of the two corners $(v,i),(u,j)$
belongs to $P$.  Row and column convexity supply the required segments.
Both corners are absent exactly when $v<L_i$ and $u>R_j$.  Such a pair exists
for these two rows exactly when $L_j<L_i$ and $R_j<R_i$.  This proves the
first assertion.  The reflected statement characterizes northeast degree
at most one.

Suppose a northwest-shifted pair $p<q$ exists, and consider a northeast-comparable pair
$a=(u,i),b=(v,j)$ which has no path with at most one bend.  Its rows are northeast shifted.
We claim that the intervals of row indices $[i,j]$ and $[p,q]$ are strictly
nested.  They cannot be disjoint: a left-endpoint increase followed by a
decrease contradicts the valley shape, and a right-endpoint decrease followed
by an increase contradicts the peak shape.  If $i<p<j<q$, the valley shape
and $L_i<L_j$ give $L_q\ge L_j$, whereas $L_p>L_q$ would make $p$ a strict
peak between $i$ and $j$.  If $p<i<q<j$, the peak shape and $R_p>R_q$ give
$R_j\le R_q$, whereas $R_i<R_j$ would make $i$ a strict valley between $p$
and $q$.  The same inequalities exclude a shared endpoint.

If $i<p<q<j$, the endpoint shapes give $L_p\le L_i$ and $R_q\ge R_j$.
Thus row $p$ contains both columns $u,v$, and the path which goes north to
row $p$, east there, and north to $b$ has two bends.  If $p<i<j<q$, they give
$R_p\le R_i$ and $L_q\ge L_j$.  Since $L_p>L_q$ and $L_p\le R_p$, we obtain
$L_j<R_i$.  The two endpoint rows therefore overlap.  Because a one-bend
path was excluded, $u<L_j$ and $v>R_i$; a common column supplies an
 east--north--east path with two bends.  Thus $D_{\mathrm{NE}}\le2$.
Reflection proves the other bound.
\end{proof}

Reflection and Lemma~\ref{pc:directional-structure} give, for every $k\ge2$,
\begin{equation}\label{eq:directional-reduction}
 Q_k=C-2(A-A_k).
\end{equation}
Indeed every excluded polyomino has exactly one directional degree above $k$,
and the other is at most one.  The two choices are disjoint and equinumerous.

\section{The core and nested column extensions}\label{pc:core-section}

A convex polyomino is \emph{directed} if it contains a cell $S$ from which
every cell is reachable by north and east steps.  This cell, necessarily its
southwest bounding-corner cell, is called the \emph{source}.  The column
containing $S$ is the \emph{source column}; its number of cells is the
\emph{source height}.

Following Beaton and Rinaldi~\cite[Section~3]{BR}, an ascending polyomino
is \emph{rectangular ascending} if it contains the northeast cell of its
bounding rectangle.  This does not require the polyomino itself to be a
rectangle.  Equivalently, the right endpoint of its top row lies in its
rightmost column.  A half-turn, that is, rotation through $180$ degrees,
identifies this class with directed convex polyominoes.
Indeed the monotone paths to the northeast corner become paths from a
southwest source.  Conversely, a convex polyomino containing its northeast
corner has nondecreasing right row endpoints by
Lemma~\ref{pc:endpoint-shapes}, hence is ascending by
Lemma~\ref{pc:directional-structure}.

\begin{definition}[Core and tail]\label{def:core-tail}
Let $P$ be ascending, with rows $0,\ldots,H_0$, where $H_0$ is the top row index.
Put
\[
 c=R_{H_0},\qquad R=P\cap\{(x,y):x\le c\}.
\]
We call $R$ the \emph{core} of $P$ and its column at $c$ the \emph{seam
column}.  This column has interval $I_0=[a_0,H_0]$.
The \emph{tail} is the ordered sequence of columns to its right,
with intervals $I_j=[a_j,b_j]$ in column $c+j$, for $1\le j\le m$.
The tail is empty when $m=0$.

Conversely, for a given rectangular ascending $R$, a \emph{nested tail}
is any finite sequence of nonempty intervals
$I_0\supseteq I_1\supseteq\cdots\supseteq I_m$ appended in successive
columns to the right of its seam column.  Such an extension is called
\emph{canonical} when the formula above recovers the given $R$ from the
resulting polyomino.  The empty extension is canonical.  No ascendingness
condition is included in the definition of a nested tail.
\end{definition}

\begin{lemma}\label{pc:core-tail}
The set $R$ is rectangular ascending, every row of $P$ meets $R$, and
$I_0\supseteq I_1\supseteq\cdots\supseteq I_m$.
Conversely, appending any nested nonempty subintervals of $I_0$ gives a convex
polyomino.
A nonempty extension is canonical exactly when its first appended column
avoids row $H_0$.  Each appended column increases the semi-perimeter by one.
\end{lemma}
\begin{proof}
A row lying wholly to the right of $c$ and the top row would be northwest
shifted.  Thus every row meets $R$.  Truncation preserves convexity and
adjacent-row overlap.  The row characterization also survives a common right
truncation, so $R$ is rectangular ascending.  Every row meeting an appended
column meets the seam; row convexity therefore forces the appended row
intervals to be initial intervals of columns.  Equivalently the $I_j$ are
nested.

Conversely nestedness makes every appended row a contiguous extension of its
old row, while each column is an interval.  The result is connected and
convex.  The top-row restriction is equivalent to
$c=R_{H_0}(P)$.  Finally a nonempty appended column leaves the height fixed and
increases the width by one.
\end{proof}

\begin{lemma}\label{pc:free-tails}
If $D_{\mathrm{NE}}(R)\ge3$, appending any nested tail gives an ascending
polyomino.
\end{lemma}
\begin{proof}
The result is convex by Lemma~\ref{pc:core-tail}.  Northeast paths between
cells of the core cannot use the appended columns, so their distances are unchanged.
Its northeast degree is therefore at least three, and
Lemma~\ref{pc:directional-structure} makes its northwest degree at most one.
\end{proof}

\begin{definition}[Shadows]\label{def:shadows}
Let $P$ be formed from a rectangular ascending core $R$ by a nested tail
as in Definition~\ref{def:core-tail}.  Its \emph{prefix of length $j$} is
the core together with the first $j$ whole tail columns.  Set $P^{[0]}=R$.  For $1\le j\le m$,
the \emph{$j$th shadow} of this extension is
\[
 P^{[j]}=P\cap\{(x,y):x\le c+j,\ y\le b_j\}.
\]
Thus it consists of the core and the first $j$ tail columns, restricted to
rows at or below the top of the $j$th column.  We call $R$ the initial shadow.
\end{definition}

Figure~\ref{fig:core-shadow} distinguishes a prefix of whole columns from its shadow and shows the half-turn of that shadow.

\begin{figure}[p]
\normalfont\normalsize
\centering
\begin{minipage}{\linewidth}
\setlength{\parindent}{0pt}
\setlength{\parskip}{5pt}
\begingroup%
\normalfont\small%
\begin{tikzpicture}[x=1cm,y=1cm,font=\small]
\path[use as bounding box] (0,.2) rectangle (15.6,17.7);

\newcommand{\pfonecore}{%
  \foreach \yy/\ll/\rr in {0/1/1,1/0/2,2/1/3,3/2/4,4/3/4,5/3/4}{%
    \foreach \xx in {\ll,...,\rr}{%
      \path[cell,fill=coreblue] (\xx,\yy) rectangle ++(1,1);}}}
\newcommand{\pfonetails}{%
  \foreach \yy in {3,4}{%
    \path[cell,fill=tailorange!65] (5,\yy) rectangle ++(1,1);
    \path[pattern=north east lines,pattern color=ink!35] (5,\yy) rectangle ++(1,1);}
  \path[cell,fill=tailgreen!65] (6,4) rectangle ++(1,1);
  \path[pattern=dots,pattern color=ink!65] (6,4) rectangle ++(1,1);}
\newcommand{\pfonecoordinates}[1]{%
  \foreach \xx in {0,...,6}{\node[font=\footnotesize] at (\xx+.5,-.4) {$\xx$};}
  \foreach \yy in {0,...,#1}{\node[font=\footnotesize] at (-.4,\yy+.5) {$\yy$};}
  \node[font=\footnotesize] at (7.35,-.4) {$x$};
  \node[font=\footnotesize] at (-.4,#1+1.35) {$y$};}

\node[anchor=west,font=\normalsize\bfseries] at (0,17.38)
  {A. Core $R$ and two tail columns};
\begin{scope}[shift={(.45,13.3)},x=.56cm,y=.56cm]
  \pfonecore
  \pfonetails
  \pfonecoordinates{5}
  \draw[ink,densely dashed,line width=.9pt] (4,3) rectangle (5,6);
  \fill[ink] (4.5,5.5) circle (1.4pt);
  \node[fill=coreblue,inner sep=.5pt,font=\footnotesize] at (4.5,3.5) {$I_0$};
  \node[fill=tailorange!65,inner sep=.5pt,font=\footnotesize] at (5.5,3.5) {$I_1$};
  \node[fill=tailgreen!65,inner sep=.5pt,font=\footnotesize] at (6.5,4.5) {$I_2$};
  \draw[decorate,decoration={brace,amplitude=2.5pt},ink]
    (4.95,-.85)--(4.05,-.85);
  \node[anchor=north,font=\footnotesize] at (4.5,-1.02) {seam $c=4$};
\end{scope}
\node[anchor=north west,text width=10cm,align=left] at (5.2,16.7) {%
  The polyomino $P$ consists of the core $R$ and the two appended columns.\\[5pt]
  $I_0=[3,5]$ in the seam column $c=4$;\\
  $I_1=[3,4]$ in column $5$;\\
  $I_2=[4,4]$ in column $6$.\\[5pt]
  The top row is $H_0=5$. The dot marks its end cell $(4,5)$; the dashed box marks the seam column.\\[5pt]
  $\operatorname{sp}(P)=7+6=13$.
};

\draw[ink!20,line width=.4pt] (0,12.25)--(15.6,12.25);

\node[anchor=west,font=\normalsize\bfseries] at (0,11.9)
  {B. The whole prefix and its second shadow};
\begin{scope}[shift={(.45,7.75)},x=.56cm,y=.56cm]
  \pfonecore
  \pfonetails
  \foreach \xx in {3,4}{%
    \path[cell,fill=ghost] (\xx,5) rectangle ++(1,1);
    \path[pattern=north west lines,pattern color=ink!35] (\xx,5) rectangle ++(1,1);}
  \pfonecoordinates{5}
  \draw[pathred,densely dashed,line width=.9pt] (-.1,5)--(7.3,5);
  \node[anchor=north,font=\footnotesize] at (3.5,-.96) {Retained cells: $P^{[2]}$};
\end{scope}
\node[anchor=north west,text width=10cm,align=left] at (5.2,11.18) {%
  The prefix of length two is all of $P$.\\[5pt]
  Truncate above row $b_2=4$ (dashed line):\\[3pt]
  $\displaystyle P^{[2]}=P\cap\{y\le4\}$.\\[5pt]
  The grey cells $(3,5)$ and $(4,5)$ are removed. Both tail columns remain in the shadow.\\[5pt]
  $Q_2=H_0-b_2=1$,\\
  $\operatorname{sp}(P^{[2]})=7+5=12$.
};

\draw[ink!20,line width=.4pt] (0,6.75)--(15.6,6.75);

\node[anchor=west,font=\normalsize\bfseries] at (0,6.4)
  {C. Half-turn of the shadow};
\begin{scope}[shift={(.45,2.52)},x=.56cm,y=.56cm]
  \foreach \yy/\ll/\rr in {0/1/1,1/0/2,2/1/3,3/2/4,4/3/4}{%
    \foreach \xx in {\ll,...,\rr}{%
      \path[cell,fill=coreblue] (6-\xx,4-\yy) rectangle ++(1,1);}}
  \foreach \yy in {3,4}{%
    \path[cell,fill=tailorange!65] (1,4-\yy) rectangle ++(1,1);
    \path[pattern=north east lines,pattern color=ink!35] (1,4-\yy) rectangle ++(1,1);}
  \path[cell,fill=tailgreen!65] (0,0) rectangle ++(1,1);
  \path[pattern=dots,pattern color=ink!65] (0,0) rectangle ++(1,1);
  \draw[ink,line width=.9pt] (0,0) rectangle (1,1);
  \node[fill=tailgreen!65,inner sep=.5pt,font=\footnotesize] at (.5,.5) {$S$};
  \pfonecoordinates{4}
  \node[anchor=north,font=\footnotesize] at (3.5,-.96) {Directed polyomino $D_2$};
\end{scope}
\node[anchor=north west,text width=10cm,align=left] at (5.2,5.76) {%
  Rotate $P^{[2]}$ through $180$ degrees and translate its source to $(0,0)$.\\[5pt]
  Cell coordinates change by\\[3pt]
  $\displaystyle (x,y)\longmapsto(6-x,4-y)$.\\[5pt]
  The last tail column becomes the source column of $D_2$: its sole cell is $S=(0,0)$, so its source height is $1$.\\[5pt]
  $\operatorname{sp}(D_2)=\operatorname{sp}(P^{[2]})=12$.
};

\path[cell,fill=coreblue] (.05,.9) rectangle +(.3,.3);
\node[anchor=west,font=\footnotesize] at (.47,1.05) {core};
\path[cell,fill=tailorange!65] (3.0,.9) rectangle +(.3,.3);
\path[pattern=north east lines,pattern color=ink!35] (3.0,.9) rectangle +(.3,.3);
\node[anchor=west,font=\footnotesize] at (3.42,1.05) {first tail column};
\path[cell,fill=tailgreen!65] (7.5,.9) rectangle +(.3,.3);
\path[pattern=dots,pattern color=ink!65] (7.5,.9) rectangle +(.3,.3);
\node[anchor=west,font=\footnotesize] at (7.92,1.05) {second tail column};
\path[cell,fill=ghost] (12.3,.9) rectangle +(.3,.3);
\path[pattern=north west lines,pattern color=ink!35] (12.3,.9) rectangle +(.3,.3);
\node[anchor=west,font=\footnotesize] at (12.72,1.05) {removed};
\end{tikzpicture}%
\endgroup%
\par
\caption{\textbf{The core, a column prefix and its shadow.} The core $R$ is blue; the first and second tail columns are distinguished by diagonal hatching and dots.
The interval $I_0$ lies in the seam column of the core.
The first appended column lowers the upper endpoint by one, and the second raises the lower endpoint by one.
The whole prefix includes the two grey cells in the top row; its shadow omits them but retains both tail columns.
Thus the shadow's semi-perimeter plus the cumulative upper truncation is $12+Q_2=13$, the semi-perimeter of the whole prefix.
The initial shadow, not truncated, is $P^{[0]}=R$. Here $D_2$ denotes the half-turn of $P^{[2]}$, translated so that its source is $(0,0)$, and $Q_2=H_0-b_2=1$ is the upper truncation. The TikZ code for this figure was prepared with assistance from ChatGPT (OpenAI) and reviewed by the author.}
\label{fig:core-shadow}
\end{minipage}
\end{figure}

Every shadow $P^{[j]}$ is convex and contains its northeast corner,
and hence is rectangular ascending.  Two cells in a shadow have exactly the
same northeast paths there as in the full polyomino: no such path can exceed
its target's row or column.  Every northeast-comparable pair lies in the shadow indexed by its
target column, or in the core.  Consequently
\begin{equation}\label{pc:shadow-maximum}
 D_{\mathrm{NE}}(P)=\max_{0\le j\le m}D_{\mathrm{NE}}(P^{[j]}).
\end{equation}
Also
\begin{equation}\label{pc:tail-degree-bound}
 D_{\mathrm{NE}}(R)\le D_{\mathrm{NE}}(P)\le D_{\mathrm{NE}}(R)+1.
\end{equation}
Old-cell distances cannot change; two northeast-comparable cells in tail columns need at most one
bend; and a cell of the core reaches a northeast-comparable cell in the tail
by first reaching the seam column in the target row and then going east.


\section{Directed polyominoes, trees, and cuts}
\label{pf:forest-section}
We use the ordered-forest representation and degree criteria of
Boussicault, Rinaldi and Socci~\cite{BRS}.  The numbered results cited below
are located in their full preprint~\cite{BRSp}.  We specify the marker
orders and cut conventions, then establish the source-row properties needed
for truncation and extension.

Translate a directed convex polyomino so that its source is $S=(0,0)$.
Its rows are intervals
$[L_y,R_y]$, numbered from the bottom.  Directedness implies that $L_y$
is nondecreasing: a path to the leftmost cell of a row must enter that row
from a column no farther east.  The distance from the source to a cell $c$ means $d_{\NE}(S,c)$.

\begin{lemma}[Distances from the source]\label{pf:source-distance}
For every integer $k\ge1$, a directed convex polyomino has degree at most
$k$ if and only if every cell can be reached from $S$ with at most $k$
bends.  In particular, if its degree is at least two, its degree is the
maximum of these source distances~\cite[Proposition~15]{BRSp}.
\end{lemma}

Write $W$ and $H$ for the width and height.  The northeast completion of a
directed convex polyomino $D$ is the polyomino $P_D$ with rows
\[
 [L_y,R_y^+],\qquad R_y^+=\max_{0\le z\le y}R_z.
\]
Its two row endpoints are nondecreasing, and consecutive rows overlap.
Such a polyomino is called a parallelogram polyomino.  Equivalently, its
two boundaries from the southwest to the northeast corner consist of north
and east steps and have no other common point.

A word is a finite ordered sequence of letters.

\begin{lemma}[Completion and cuts]\label{pf:completion-cut}
Completion changes no source-to-cell distance for cells of $D$.
If $a$ is the left endpoint of the top row of a parallelogram $P$, and
$b$ is the bottom endpoint of its rightmost column, directed polyominoes
with completion $P$ are in bijection with words having $W-a$ east steps
and $H-b$ south steps, starting east and ending south.
\end{lemma}
\begin{proof}
We first show that no newly filled cell is southwest of an old cell.
Suppose $(x,y)$ is filled because $x>R_y$ but $x\le R_z$ for some $z<y$.
Since $L_z\le L_y\le R_y<x$, column $x$ contains a cell below row $y$.
If an old cell $c$ had $c_x\ge x$ and $c_y\ge y$, a directed path to
$c$ would cross row $y$ in a column at most $R_y<x$ and subsequently
cross column $x$ above row $y$.  This would contradict column convexity.
Consequently no north-east path to an old target can use a filled cell.
All its source distances are unchanged.

The cut bijection is Proposition~3 of~\cite{BRSp}, with its
northeast corner rectangle written as $[a,W-1]\times[b,H-1]$ in cell
coordinates.  The full cut runs from $(a,H)$ to $(W,b)$; its first east
and last south steps retain the top-left and bottom-right cells of this
rectangle.  These are the endpoint conventions used below.
Inside the rectangle, $D$ consists of the cells southwest of the cut;
outside it, $D$ agrees with $P$.
\end{proof}

A \emph{rooted tree} is a finite connected graph without cycles, with a
distinguished vertex called its root.  The \emph{parent} of a non-root
vertex is the next vertex on its path to the root; its other neighbours
are its \emph{children}; all neighbours of the root are its children.
A \emph{leaf} has no children.  An
\emph{ordered rooted tree} has a linear order on the children of each
vertex.  The \emph{depth} $\operatorname{dep}(v)$ counts the vertices on
the path from $v$ to the root, so a root has depth one.  Vertices of depth
$j$ form \emph{level $j$}, and the greatest depth is the \emph{height}.
The path from a vertex to its root, listed in that order, is its
\emph{ancestor chain}; its length also counts vertices.

An \emph{ordered forest} is a finite ordered list of these trees, including
the empty list.  Its height is the greatest tree height, or zero if empty.
Within each tree, level one consists of the root alone.  At each later level, order the vertices first by the order of their
parents at the preceding level and then by the specified order of children
of one parent.  This recursively defines the \emph{level order}.

We now represent a parallelogram by two ordered forests.  Its column markers are the top cells of all columns
except the rightmost; its row markers are the rightmost cells of all rows
except the top.  The northeast cell $E=(W-1,H-1)$ is not marked.
These two marker sets are disjoint.  Indeed, the rightmost cell of a non-top
row has a cell immediately above it, and the top cell of a non-rightmost
column has a cell immediately to its right.

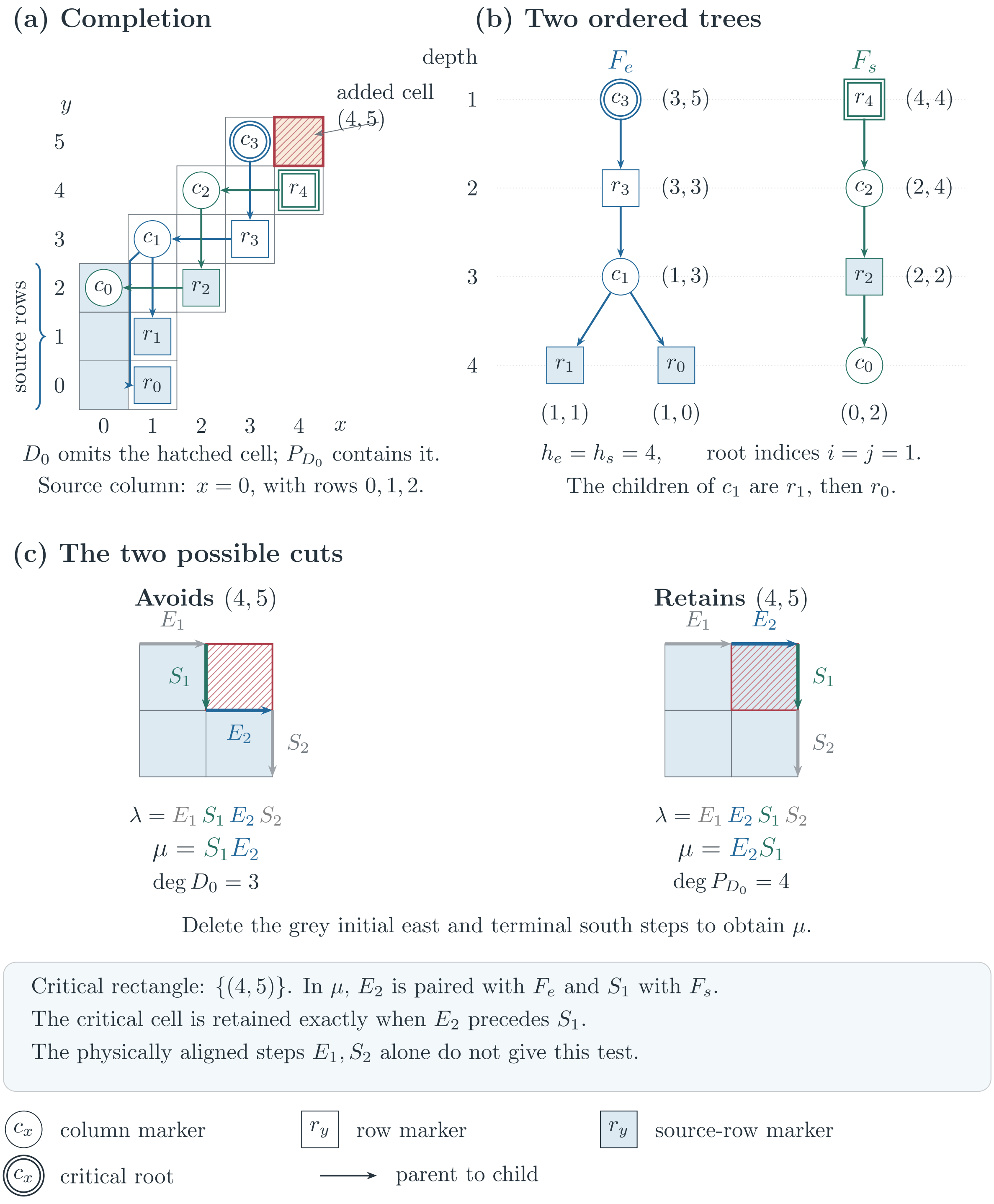
\begin{figure}[p]
\normalfont\normalsize
\centering
\begin{minipage}{\linewidth}
\setlength{\parindent}{0pt}
\setlength{\parskip}{5pt}
\begin{tikzpicture}[x=1cm,y=1cm,
  colmark/.style={circle,draw=ink,fill=white,minimum size=5.8mm,inner sep=.4pt,font=\small},
  rowmark/.style={rectangle,draw=ink,fill=white,minimum size=5.8mm,inner sep=.4pt,font=\small},
  criticalmark/.style={line width=.7pt,double,double distance=.9pt},
  sourcemark/.style={fill=coreblue},
  eastedge/.style={draw=wallblue,line width=.9pt,-{Stealth[length=1.5mm]}},
  southedge/.style={draw=safeolive,line width=.9pt,-{Stealth[length=1.5mm]}},
  coord/.style={font=\small,anchor=north},
  cutend/.style={draw=ink!45,line width=1.4pt,-{Stealth[length=1.8mm]}},
  cute/.style={draw=wallblue,line width=1.5pt,-{Stealth[length=1.8mm]}},
  cuts/.style={draw=safeolive,line width=1.5pt,-{Stealth[length=1.8mm]}}]
\path[use as bounding box] (0,-.28) rectangle (15.6,18.65);
\node[font=\bfseries\large,anchor=west] at (0,18.3) {(a) Completion};
\node[font=\bfseries\large,anchor=west] at (7.3,18.3) {(b) Two ordered trees};
\begin{scope}[yshift=6.8cm]
\begin{scope}[shift={(1.2,5.35)},scale=.77]
  \foreach \y/\a/\b in {0/0/1,1/0/1,2/0/2,3/1/3,4/2/4,5/3/4}
    \foreach \x in {\a,...,\b}
      \path[cell,fill=white] (\x,\y) rectangle +(1,1);
  \foreach \y in {0,1,2}
    \path[cell,fill=coreblue] (0,\y) rectangle +(1,1);
  \path[cell,fill=tailorange!22] (4,5) rectangle +(1,1);
  \path[pattern=north east lines,pattern color=pathred!65] (4,5) rectangle +(1,1);
  \draw[pathred,line width=1.15pt] (4,5) rectangle +(1,1);
  \foreach \x in {0,...,4} \node[font=\small] at (\x+.5,-.33) {$\x$};
  \foreach \y in {0,...,5} \node[font=\small,anchor=east] at (-.12,\y+.5) {$\y$};
  \node[font=\small,anchor=west] at (5.05,-.33) {$x$};
  \node[font=\small] at (-.27,6.2) {$y$};
  \draw[decorate,decoration={brace,mirror,amplitude=4pt},wallblue,line width=.8pt]
    (-.9,0)--(-.9,3);
  \node[font=\small,rotate=90] at (-1.2,1.5) {source rows};
  \node[colmark,draw=wallblue,criticalmark] (gc3) at (3.5,5.5) {$c_3$};
  \node[rowmark,draw=wallblue] (gr3) at (3.5,3.5) {$r_3$};
  \node[colmark,draw=wallblue] (gc1) at (1.5,3.5) {$c_1$};
  \node[rowmark,draw=wallblue,sourcemark] (gr1) at (1.5,1.5) {$r_1$};
  \node[rowmark,draw=wallblue,sourcemark] (gr0) at (1.5,.5) {$r_0$};
  \node[rowmark,draw=safeolive,criticalmark] (gr4) at (4.5,4.5) {$r_4$};
  \node[colmark,draw=safeolive] (gc2) at (2.5,4.5) {$c_2$};
  \node[rowmark,draw=safeolive,sourcemark] (gr2) at (2.5,2.5) {$r_2$};
  \node[colmark,draw=safeolive] (gc0) at (.5,2.5) {$c_0$};
    \draw[eastedge] (gc3)--(gr3);
    \draw[eastedge] (gr3)--(gc1);
    \draw[eastedge] (gc1)--(gr1);
    \draw[eastedge] (gc1.south west) -- (1.04,3.05) -- (1.04,.5) -- (gr0.west);
    \draw[southedge] (gr4)--(gc2);
    \draw[southedge] (gc2)--(gr2);
    \draw[southedge] (gr2)--(gc0);
  \draw[ink!65,-{Stealth[length=1.5mm]}] (6.15,5.9) -- (4.82,5.64);
  \node[font=\small,anchor=west,align=left] at (5.1,6.2)
    {added cell\\$(4,5)$};
\end{scope}
\node[font=\small,align=center] at (3.6,4.62)
  {$D_0$ omits the hatched cell; $P_{D_0}$ contains it.};
\node[font=\small,align=center] at (3.6,4.13)
  {Source column: $x=0$, with rows $0,1,2$.};

\end{scope}
\begin{scope}[xshift=-.4cm,yshift=6.8cm]
\node[font=\large\bfseries,text=wallblue] at (10.15,10.85) {$F_e$};
\node[font=\large\bfseries,text=safeolive] at (14.0,10.85) {$F_s$};
\node[font=\small,anchor=east] at (8.03,10.9) {depth};
\foreach \d/\y in {1/10.25,2/8.85,3/7.45,4/6.05} {
  \draw[ink!15,densely dotted] (8.2,\y)--(15.6,\y);
  \node[font=\small,anchor=east] at (8.03,\y) {$\d$};
}
\node[colmark,draw=wallblue,criticalmark] (ec3) at (10.15,10.25) {$c_3$};
\node[rowmark,draw=wallblue] (er3) at (10.15,8.85) {$r_3$};
\node[colmark,draw=wallblue] (ec1) at (10.15,7.45) {$c_1$};
\node[rowmark,draw=wallblue,sourcemark] (er1) at (9.27,6.05) {$r_1$};
\node[rowmark,draw=wallblue,sourcemark] (er0) at (11.03,6.05) {$r_0$};
\node[rowmark,draw=safeolive,criticalmark] (sr4) at (14,10.25) {$r_4$};
\node[colmark,draw=safeolive] (sc2) at (14,8.85) {$c_2$};
\node[rowmark,draw=safeolive,sourcemark] (sr2) at (14,7.45) {$r_2$};
\node[colmark,draw=safeolive] (sc0) at (14,6.05) {$c_0$};
\draw[eastedge] (ec3)--(er3);
\draw[eastedge] (er3)--(ec1);
\draw[eastedge] (ec1)--(er1);
\draw[eastedge] (ec1)--(er0);
\draw[southedge] (sr4)--(sc2);
\draw[southedge] (sc2)--(sr2);
\draw[southedge] (sr2)--(sc0);
\foreach \v/\x/\y in {ec3/3/5,er3/3/3,ec1/1/3,sr4/4/4,sc2/2/4,sr2/2/2}
  \node[font=\small,anchor=west] at ([xshift=2mm]\v.east) {$(\x,\y)$};
\foreach \v/\x/\y in {er1/1/1,er0/1/0,sc0/0/2}
  \node[coord] at ([yshift=-1.6mm]\v.south) {$(\x,\y)$};
\node[font=\small] at (11.9,4.65) {$h_e=h_s=4$,\qquad root indices $i=j=1$.};
\node[font=\small] at (11.9,4.12) {The children of $c_1$ are $r_1$, then $r_0$.};

\end{scope}
\node[font=\bfseries\large,anchor=west] at (0,9.85) {(c) The two possible cuts};
\node[font=\bfseries] at (3.2,9.15) {Avoids $(4,5)$};
\node[font=\bfseries] at (11.5,9.15) {Retains $(4,5)$};
\foreach \xx in {2.15,10.45} {
  \begin{scope}[shift={(\xx,6.35)},scale=1.05]
    \foreach \x in {0,1} \foreach \y in {0,1}
      \path[cell,fill=coreblue] (\x,\y) rectangle +(1,1);
  \end{scope}
}
\begin{scope}[shift={(2.15,6.35)},scale=1.05]
  \path[cell,fill=white] (1,1) rectangle +(1,1);
  \path[pattern=north east lines,pattern color=pathred!65] (1,1) rectangle +(1,1);
  \draw[pathred,line width=.8pt] (1,1) rectangle +(1,1);
  \draw[cutend] (0,2)--(1,2) node[midway,above=2pt,font=\small,text=ink!65] {$E_1$};
  \draw[cuts] (1,2)--(1,1) node[midway,left=2pt,font=\small,text=safeolive] {$S_1$};
  \draw[cute] (1,1)--(2,1) node[midway,below=2pt,font=\small,text=wallblue] {$E_2$};
  \draw[cutend] (2,1)--(2,0) node[midway,right=2pt,font=\small,text=ink!65] {$S_2$};
\end{scope}
\begin{scope}[shift={(10.45,6.35)},scale=1.05]
  \path[pattern=north east lines,pattern color=pathred!65] (1,1) rectangle +(1,1);
  \draw[pathred,line width=.8pt] (1,1) rectangle +(1,1);
  \draw[cutend] (0,2)--(1,2) node[midway,above=2pt,font=\small,text=ink!65] {$E_1$};
  \draw[cute] (1,2)--(2,2) node[midway,above=2pt,font=\small,text=wallblue] {$E_2$};
  \draw[cuts] (2,2)--(2,1) node[midway,right=2pt,font=\small,text=safeolive] {$S_1$};
  \draw[cutend] (2,1)--(2,0) node[midway,right=2pt,font=\small,text=ink!65] {$S_2$};
\end{scope}
\node[font=\small] at (3.2,5.73)
  {$\lambda={\color{gray}E_1}\,{\color{safeolive}S_1}\,{\color{wallblue}E_2}\,{\color{gray}S_2}$};
\node[font=\small] at (11.5,5.73)
  {$\lambda={\color{gray}E_1}\,{\color{wallblue}E_2}\,{\color{safeolive}S_1}\,{\color{gray}S_2}$};
\node[font=\large] at (3.2,5.2) {$\mu={\color{safeolive}S_1}{\color{wallblue}E_2}$};
\node[font=\large] at (11.5,5.2) {$\mu={\color{wallblue}E_2}{\color{safeolive}S_1}$};
\node[font=\small] at (3.2,4.67) {$\deg D_0=3$};
\node[font=\small] at (11.5,4.67) {$\deg P_{D_0}=4$};
\node[font=\small,align=center] at (7.8,3.98)
  {Delete the grey initial east and terminal south steps to obtain $\mu$.};

\draw[rounded corners=2mm,fill=coreblue!35,draw=ink!25] (0,1.38) rectangle (15.6,3.42);
\node[font=\small,anchor=west] at (.3,3.02)
  {Critical rectangle: $\{(4,5)\}$. In $\mu$, $E_2$ is paired with $F_e$ and $S_1$ with $F_s$.};
\node[font=\small,anchor=west] at (.3,2.5)
  {The critical cell is retained exactly when $E_2$ precedes $S_1$.};
\node[font=\small,anchor=west] at (.3,1.98)
  {The physically aligned steps $E_1,S_2$ alone do not give this test.};
\node[colmark] at (.32,.78) {$c_x$};
\node[font=\small,anchor=west] at (.76,.78) {column marker};
\node[rowmark] at (5.0,.78) {$r_y$};
\node[font=\small,anchor=west] at (5.44,.78) {row marker};
\node[rowmark,sourcemark] at (9.72,.78) {$r_y$};
\node[font=\small,anchor=west] at (10.16,.78) {source-row marker};
\node[colmark,criticalmark] at (.32,.05) {$c_x$};
\node[font=\small,anchor=west] at (.76,.05) {critical root};
\draw[ink,line width=.8pt,-{Stealth[length=1.6mm]}] (5.0,.05)--(5.9,.05);
\node[font=\small,anchor=west] at (6.06,.05) {parent to child};
\end{tikzpicture}
\par
\caption{\textbf{Completion, forests and the cut.} Let $D_0$ be the half-turn of the core $R$ in Figure~\ref{fig:core-shadow}, translated to put its source at $(0,0)$. Its rows, from bottom to top, are
$[0,1],[0,1],[0,2],[1,3],[2,4],[3,3]$.
Its completion adds one cell, which is also the entire critical rectangle.
The top-right cell has no marker.  Root depth is one, and the two trees have height four.
In (c), blue cells are retained.  The cut comparison uses the letters paired with the trees in $\mu$, not the steps physically aligned with their roots. Shaded row markers belong to rows meeting the source column; they are called source-row markers. The TikZ code for this figure was prepared with assistance from ChatGPT (OpenAI) and reviewed by the author.}
\label{fig:completion-cut}
\end{minipage}
\end{figure}

A marker in the top row is a root of the forest $F_e$; a marker in the
rightmost column is a root of $F_s$.  The parent of every other column
marker is the row marker in its row, and the parent of every other row
marker is the column marker in its column.  Parents lie strictly north or
east, so these relations form trees.  Children in a row are ordered from
east to west, and children in a column from north to south.  The roots of
$F_e$ are ordered west to east, and those of $F_s$ north to south.
Consequently column markers are the odd-depth
vertices of $F_e$ and the even-depth vertices of $F_s$; row markers have the
opposite parities.
The marker construction is illustrated in Figure~\ref{fig:completion-cut}(a,b).

For a parallelogram $P$, let $R(y)$ be the right endpoint of row $y$
and $T(x)$ the top endpoint of column $x$.  Write
\[
 r_y=(R(y),y)\quad(0\le y<H-1),\qquad
 c_x=(x,T(x))\quad(0\le x<W-1)
\]
for its row and column markers.  All forest terms for a directed polyomino
$D$ refer to the forests of $P_D$; in particular $R(y)=R_y^+$ there.

\begin{lemma}[Order of marker levels]\label{pf:row-order}
From bottom to top, the depths of the row markers are nonincreasing.
At a fixed depth their roots occur in increasing coordinate order:
west to east for $F_e$, and south to north for $F_s$.
Within a single tree, its row markers at a fixed depth occur from top to
bottom in the level order.

The corresponding assertions for column markers hold from left to right:
their depths are nonincreasing; at one depth their roots occur west to east
in $F_e$ or south to north in $F_s$; and within a tree their level order
is east to west.
\end{lemma}
\begin{proof}
The functions $R$ and $T$ are nondecreasing.  Away from the roots the
parent maps alternate as
\[
 r_y\longmapsto c_{R(y)},\qquad
 c_x\longmapsto r_{T(x)}.
\]
Start two chains at rows $y<z$.  Applying $R$ and $T$ preserves the order
of their coordinates at every stage.  If the lower chain reaches the right
boundary, the upper one has done so; if it reaches the top boundary, the
upper one has done so.  The upper chain therefore has no greater depth.
If the depths are equal, the same coordinate comparison persists to the
roots and proves their stated order.  The argument starting at two columns
is identical with the roles of $R$ and $T$ exchanged.

For the within-tree assertion, read a column-marker level from east to
west.  Its row children are grouped by their parents; monotonicity of $R$
places all children of the more eastern parent above those of the more
western parent.  Within a group they are already ordered north to south.
Thus the next row level is read top to bottom.  From a row level, the same
argument with $T$ gives east-to-west order at the next column level.
Induction from the root proves the claim.
\end{proof}

At a fixed depth, the row markers belonging to one tree form a consecutive
block in vertical order.  We call it that tree's \emph{row-marker block}
at the specified depth.

\begin{proposition}[The forest bijection]\label{pf:forest-bijection}
The preceding construction is a bijection from parallelogram polyominoes
to pairs of ordered forests.  A pair with $n$ vertices corresponds to
semi-perimeter $n+2$.  The numbers of trees in $F_e,F_s$ are the top-row
length minus one and the rightmost-column length minus one, respectively~\cite[Proposition~2]{BRSp}.
\end{proposition}

For later degree calculations, use two maximal north-east paths from $S$
to a target, one starting east and the other north.  A run is a maximal sequence of consecutive steps in one direction.
At each run these paths go as far as possible within the polyomino toward
the target, without passing its row or column.  Include an incoming unit segment
before $S$ in the chosen initial direction.  This convention allows a
turn at $S$ when the initial internal run has length zero.  It does not
change the smaller bend count of the two paths: an ordinary path can be
given an incoming segment in its actual initial direction without creating
a bend, and deleting an incoming segment cannot increase its bend count.

These paths are also called \emph{bounce paths}.  For target $E$, they
are the two \emph{maximal boundary paths}.

\begin{lemma}[Maximal paths]\label{pf:maximal-paths}
The smaller bend count of these two maximal paths is the minimum bend
count between their endpoints.  In a parallelogram the largest source
distance is the distance from $S$ to $E$~\cite[Lemmas~16, 18 and 19]{BRSp}.
\end{lemma}

Let the height of a forest be the greatest vertex depth, with height zero
for the empty forest.  Write these heights as $h_e,h_s$.

The degree formula below is Proposition~22 of~\cite{BRSp}.
The following derivation also identifies the extremal roots
under our marker order.
\begin{proposition}[Degree and the two forest heights]\label{pf:bounce-degree}
The degree of a parallelogram is
\[
 \max(h_e,h_s)-\begin{cases}0,&h_e=h_s,\\1,&h_e\ne h_s.\end{cases}
\]
When $h_e=h_s=h>0$, both maximal boundary paths have $h$ bends.
Their terminal roots are the roots of the first height-$h$ tree of $F_e$ and the last
height-$h$ tree of $F_s$.  We call this the \emph{flat case}.
If the heights are unequal with maximum $h$, their bend counts are
$h$ and $h-1$, and both chains lie in the taller forest.
\end{proposition}
\begin{proof}
If $W=1$ or $H=1$, the polyomino has a single row or column (a \emph{bar}), and the
formula follows immediately from its empty forests or depth-one roots.
Assume henceforth that $W,H\ge2$, so both extreme markers exist.
The north-first path turns first at the source-column top and thereafter
at its successive parents.  Its bend count is therefore the depth of the
leftmost column marker.  Similarly the east-first count is the depth of
the bottom row marker.  These assertions also hold for an internal initial
run of length zero, by the incoming-segment convention.

Lemma~\ref{pf:row-order} says that these depths are the largest column-
and row-marker depths.  If both forest heights are $h$, the parity rule
places depth-$h$ column markers in one forest and depth-$h$ row markers
in the other.  Both maxima are $h$.  Their roots are, by the same lemma,
the roots of the first height-$h$ east tree and the last height-$h$ south tree.

If only one forest has height $h$, its vertices at depths $h$ and $h-1$
give one marker maximum equal to $h$ and the other equal to $h-1$.
The other forest cannot increase the second maximum: its depth-$h-1$
vertices, if present, have the first marker type.  Thus the counts are
$h,h-1$.  Lemmas~\ref{pf:source-distance} and~\ref{pf:maximal-paths}
give the formula, including a rectangle when the common height is one.
\end{proof}

A parallelogram is \emph{flat} when $h_e=h_s>0$, and its \emph{flat
height} is this common forest height.  A directed convex polyomino is flat
when its completion is flat.  The two terminal roots in the flat case are
its \emph{critical east and south roots}; their trees are the
\emph{critical trees}.  Otherwise we use the term \emph{nonflat}.

In the flat case, the penultimate intersection of the two maximal paths
is $(c,r)$ when their terminal roots are $(c,H-1)$ and $(W-1,r)$.
This identifies the coordinates in the rectangle criterion of~\cite{BRSp}.
\begin{lemma}[The highest-degree rectangle]\label{pf:critical-rectangle}
In the flat case of height $h\ge2$, let the two terminal roots have
coordinates $(c,H-1)$ and $(W-1,r)$.  A cell $(x,y)$ of the parallelogram
has source distance $h$ exactly when
\[
 x>c\quad\hbox{and}\quad y>r.
\]
Every other cell has source distance at most $h-1$~\cite[Corollary~29 and Lemma~30]{BRSp}.
\end{lemma}

The cells satisfying these two strict inequalities form the
\emph{critical rectangle}.  A cut \emph{meets} this rectangle when its
retained polyomino contains a cell of the rectangle; otherwise it
\emph{avoids} it.

\begin{proposition}[Degree after a cut]\label{pf:completion-degree}
Let $D$ have completion $P$ of degree $q\ge2$.
If $P$ is not flat then $D$ has degree $q$.
If $P$ is flat then $D$ has degree $q$ when it meets the rectangle of
Lemma~\ref{pf:critical-rectangle}, and degree $q-1$ otherwise~\cite[Proposition~31]{BRSp}.
\end{proposition}
For $q=2$, an avoiding cut has degree at most one by the cited result and
at least one because its unchanged width and height are both greater than
one.  Thus the exact-degree conclusion also covers this endpoint.

Let $\lambda$ be the full cut, and let $m_e,m_s$ be the numbers of trees
in $F_e,F_s$.  Number its east steps $E_1,\ldots,E_{m_e+1}$ and south steps
$S_1,\ldots,S_{m_s+1}$ in their respective orders.  Delete its first east
and last south steps and write
\[
 \lambda=E_1\,\mu\,S_{m_s+1}.
\]
We call $\mu$ the \emph{cut word}.  It is an arbitrary sequence of
$m_e$ east letters and $m_s$ south letters.  Pair its $i$th east letter with the $i$th tree of $F_e$, and
its $j$th south letter with the $j$th tree of $F_s$.  Thus their actual
positions in the full cut are $E_{i+1}$ and $S_j$.  This pairing is used
only for tree order and the following rectangle test.  The steps physically
in the critical root column and row would instead be $E_i,S_{j+1}$.

\begin{lemma}[Critical order in the cut word]\label{pf:critical-order}
Suppose $P$ has flat height $h\ge2$.  If its critical east and south trees have root indices
$i,j$, then $D$ meets the highest-degree rectangle if and only if the
$i$th east letter precedes the $j$th south letter in $\mu$.
\end{lemma}
\begin{proof}
The top-row roots occupy columns $a,a+1,\ldots,W-2$, and the right-column
roots occupy rows $H-2,H-3,\ldots,b$.  Hence the critical coordinates are
$c=a+i-1$ and $r=H-1-j$.  The part of $D$ in the northeast corner
rectangle is southwest closed.  It meets $x>c,y>r$ precisely when it
contains $(c+1,r+1)$.  The east step above this column is $E_{i+1}$.
That cell is retained precisely when this east step is taken before the
cut descends to height $r+1$, namely before $S_j$.  These are exactly
the two paired letters in $\mu$.
\end{proof}

Figure~\ref{fig:completion-cut} illustrates the criterion for the half-turn $D_0$ of the core $R$ in Figure~\ref{fig:core-shadow}. The two cut words $S_1E_2$ and $E_2S_1$ respectively omit and retain the single cell of the critical rectangle.

The prefixes of a word are its
initial sequences, and its suffixes its terminal sequences; both may be
empty.  In the cut words used for enumeration, each letter carries one tree.
A tree with $n$ vertices has weight $x^n$.  To mark a nonnegative
integer-valued statistic $s$, multiply the weight by $u^s$.  Word
weights are products of letter weights, and generating functions are sums
of these weights.  We use the same letter for a class and its generating
function when the context distinguishes them.

The bijections above give cut generating functions with arbitrary tree weights.  If east and south trees have weights $e,s$, then
\[
 \frac{x^2}{1-e-s}
\]
counts the forest pairs and cuts; $x^2$ supplies the two further units of semi-perimeter in
Proposition~\ref{pf:forest-bijection}.
If one specified tree has weight $M$, with tree weights $e_-,s_-$ before
its letter and $e_+,s_+$ after it, its series is
\begin{equation}\label{pf:marked-cut-series}
 \frac{x^2M}{(1-e_--s_-)(1-e_+-s_+)}.
\end{equation}
These statements follow simply by expanding each denominator as a sum of
finite words.  They remain valid when each tree weight records additional
statistics of that tree.

For a set of letter types $\mathcal A$, the notation $\mathcal A^*$
denotes all finite words in these types, including the empty word.
Juxtaposition of word classes denotes concatenation in the stated order.

For a parallelogram of flat height $h\ge2$, let $e_1,s_1$ designate trees of
height at most $h-1$, and $e_2,s_2$ trees of height at most $h$.
Let $u_e,u_s$ designate the two distinguished letters carrying the
critical trees, each of height $h$.  Avoidance of the highest-degree
rectangle is equivalent to the following decomposition of $\mu$:
\begin{equation}\label{pf:flat-avoiding-word}
 \{e_1,s_2\}^*u_s\{e_1,s_1\}^*u_e\{e_2,s_1\}^*.
\end{equation}
Indeed, Lemma~\ref{pf:critical-order} requires $u_s$ before $u_e$;
within-direction root order gives the bounds in each region.  Consequently
its fully weighted generating function is
\begin{equation}\label{pf:flat-avoiding-series}
 \frac{x^2u_eu_s}
 {(1-e_1-s_2)(1-e_1-s_1)(1-e_2-s_1)}.
\end{equation}
Weights and marks on specified trees or levels can therefore be inserted
in these regions without changing the cut interpretation.

Let $D$ have source height $r$.  Its \emph{source rows} are the rows
$0\le y<r$; completion leaves this range unchanged.  When $r<H$, these
rows have markers in $P_D$, called the \emph{source-row markers}; the depth
of a source row means the depth of its marker.  In this case the highest cell of its source column is
$w=c_0=(0,r-1)$, called the \emph{source top}.  The tree containing $w$,
its root, and its forest are respectively the \emph{source tree},
\emph{source root} and \emph{source forest}.  The parent $r_{r-1}$ of $w$
is the \emph{source parent}.

\begin{lemma}[The source tree and its source rows]\label{pf:source-extremal}
In the completion $P_D$, assume the source column ends below the top
row, at a marker $w$ of depth
$d\ge2$.  Its parent is the bottommost non-leaf row marker at depth
$d-1$, and is the last non-leaf at that level in its tree's level order.
The vertex $w$ is its last child.  The source tree is the first tree of
$F_e$ reaching depth $d$, if $w\in F_e$, or the last such tree of $F_s$
if $w\in F_s$.

All row-marker levels deeper than $d-1$ lie entirely in the source rows.  At level $d-1$, all blocks below the source tree's block are wholly
in the source rows, all blocks above it are outside those rows, and the retained markers
in the source tree are the parent of $w$ and the later vertices in
that level order.
\end{lemma}
\begin{proof}
Below $w$, a row and the next row both start in column zero, while their
right endpoints are nondecreasing.  Every cell in the lower row therefore
has a cell above it.  There is no column marker in such a row, so its row
marker has no children.  The row of $w$ does have a column marker, and
its row marker is the parent of $w$.  It is thus the bottommost non-leaf
at its level.  The child order in its row makes the leftmost column marker
$w$ its last child.  Lemma~\ref{pf:row-order} turns bottommost into last
in the within-tree level order.  The column-marker part of that lemma
places the leftmost marker $w$ at the greatest column-marker depth and
at the first applicable root in coordinate order.  Trees of its forest
reach depth $d$ exactly when they have a column marker at that depth,
which proves the stated extremal-root assertion.

The same lemma makes each row-marker depth a consecutive block.  All
deeper blocks lie below the parent's row and hence inside the bottom-anchored
source rows.  At its own depth, row-marker blocks are consecutive and only the
source tree's block can meet both the source rows and the rows above them.  Its within-tree order
then gives exactly the retained markers stated above.
\end{proof}

In a nonflat degree-$q$ parallelogram with $q\ge2$ the source-top depth is $q$ or
$q+1$; in a parallelogram of flat height $h\ge2$ it is $h$.
These assertions follow directly from the two boundary-chain lengths in
Proposition~\ref{pf:bounce-degree}, since the source top starts the
north-first chain.  Together with Lemma~\ref{pf:source-extremal} and the
weighted word formulas, they provide the source-column refinements used
below.

\begin{lemma}[Adding a source column]\label{pf:source-extension}
Let $D$ be directed convex with source height $r<H$, and adjoin a new leftmost column of
height $1\le n\le r$, anchored at the bottom.  The new completion is the
old completion translated one unit east, together with the new column.
Its forest is obtained by adjoining one leaf: the new column marker
$(0,n-1)$ is the last child of the translated row marker in row $n-1$.
Every old marker, parent and root order is unchanged, as is the cut word.
The new north-first boundary path has bend count
$1+\operatorname{dep}(r_{n-1})$.
\end{lemma}
\begin{proof}
Every old row endpoint is translated one unit east, so its running right
maximum is translated by the same amount.  The completion formula therefore
gives precisely the asserted new parallelogram.  Each new cell has an old
cell immediately to its right.  Thus no new cell is a row marker, and the
only new column marker is $(0,n-1)$; all old topmost and rightmost cells are
unchanged.  Since $n<H$, this new marker is not a root.  Its parent is the
old row marker in its row.  No row marker lies in its column, so it is a
leaf, and the east-to-west child order makes it the last child.  The
northeast cut rectangle and the retained cells in it are merely translated,
so its word is unchanged.  The north-first path follows the new marker and
its ancestor chain; the incoming-segment convention also covers $n=1$.
\end{proof}

\begin{lemma}[Deleting bottom rows]\label{pf:bottom-pruning}
Delete $q<r$ bottom rows of a directed convex polyomino of source height $r$.
If the old rightmost column survives, completion commutes with this deletion,
up to translation.  In the forest this deletes exactly the row markers in
the removed rows, and every retained marker keeps its parent.  If none of
the deleted row markers is a root, the root orders and cut word are unchanged.
\end{lemma}
\begin{proof}
Write $M=\max_{y<q}R_y$ when $q>0$.  A retained cell of the rightmost
column exists.  If $R_q<M$, a directed path to that cell crosses row $q$
strictly west of column $M$ and later crosses column $M$ above that row.
That column also meets a deleted row, contrary to column convexity.
Hence $R_q\ge M$, and taking running right maxima commutes with deleting
the first $q$ terms.  This proves the completion assertion.

In the completed polyomino, every row below $q$ and its successor start
at column zero, and their right endpoints are nondecreasing.  No such row
contains a column marker.  Thus the deletion removes only row markers,
all leaves, and every remaining column still has its original top cell.
Every parent lies in the same row or strictly above its child, so it
survives with the same coordinates.  Finally, a deletion which removes no
root lies below the bottom of the rightmost column.  It therefore misses
the northeast cut rectangle and preserves the cut word.
\end{proof}

\section{Cumulative truncation and the first forbidden column}\label{sec:truncation}

For a directed convex polyomino $D$ of source height $r$, and integers
$0\le q<r$, let $\operatorname{Trim}_q(D)$ delete the bottom $q$ rows
and translate the remainder down by $q$.  For $0\le p<r$, let
$\operatorname{Ext}_p(D)$ translate $D$ one column east and adjoin the
column $\{(0,y):0\le y<r-p\}$.

For a nested tail with $I_j=[a_j,b_j]$ and $b_0=H_0$, let $D_j$ be its
$j$th shadow after a half-turn and translation to source $(0,0)$.  Its
source height is $r_j=b_j-a_j+1$, and put $r=r_0$.  The endpoint changes
\[
 p_j=a_{j+1}-a_j,\qquad q_j=b_j-b_{j+1}
\]
are nonnegative integers.  They give
\[
 D_{j+1}=\operatorname{Ext}_{p_j}(\operatorname{Trim}_{q_j}(D_j)),
 \qquad p_j+q_j<r_j.
\]
We call any finite sequence of such choices an \emph{extension sequence},
including the empty sequence.  It is canonical when it is empty or its
first choice satisfies $q_0\ge1$.  The \emph{cumulative slacks} are
\[
 P_j=\sum_{i<j}p_i=a_j-a_0,\qquad
 Q_j=\sum_{i<j}q_i=H_0-b_j,
\]
so $r_j=r-P_j-Q_j$.  For a fixed bound $k$, a sequence has its
\emph{first passage} at an index $j\ge1$ when $\deg D_i\le k$ for $0\le i<j$ and
$\deg D_j>k$.  Its $j$th shadow is then called the first forbidden shadow.

\begin{lemma}\label{pc:history-contraction}
For a nonempty extension sequence,
\begin{equation}\label{eq:history-contraction}
 \deg D_j=\deg\operatorname{Ext}_{P_j}(\operatorname{Trim}_{Q_j}(D_0)).
\end{equation}
The empty extension sequence remains $D_0$.  Duplicating a source column preserves degree
whenever the old degree is at least one.
\end{lemma}
\begin{proof}
After the cumulative bottom trim, the earlier inserted columns are
bottom-anchored and their heights increase towards the old core.  Each
contains every row of the newest source column.  Contracting the intermediate
columns leaves the right side of \eqref{eq:history-contraction}.
A path in the contracted polyomino expands by prolonging its horizontal passage
through these columns.  Conversely let the newest source column have height $n$ and let
$(j,y)$ be the first cell at which an old path enters the old core; its source
column is at abscissa $j$.  If $y<n$, replace the prefix by $N^yE^j$.
This ends east, like the old prefix, and uses at most one bend, no more than
the old prefix.  If $y\ge n$, the old prefix must leave the newest column
eastwards before rising, and must again move east to enter the old core;
it has at least two bends.  Replace it by $E^jN^y$ in the bottom row and
old source column.  The latter contains row $y$, and therefore the whole
vertical segment.  This replacement has one bend and at most one joining
bend, again no more than the old prefix.  Thus source distances to old-core
cells agree.  An intermediate-column target needs at
most one bend, already attained in the contracted polyomino unless it is a bar.
The bar cases are immediate.  The source-distance characterization in
Lemma~\ref{pf:source-distance} gives the equality of degrees.

For a duplicated column the same argument applies to targets outside the old
source column; targets inside it need at most one bend.  Compression gives
the reverse inequality.  The qualification on degree excludes a vertical bar
of height greater than one, whose duplication has degree one.
\end{proof}

\begin{lemma}\label{pc:slack-monotonicity}
Put $d=\deg D_0$ and
\[
 \delta(P,Q)=\deg\operatorname{Ext}_P(\operatorname{Trim}_Q(D_0)),
 \qquad P,Q\ge0,\quad P+Q<r.
\]
Then $\delta(P,Q)\in\{d-1,d,d+1\}$; it is nondecreasing in $P$ and
nonincreasing in $Q$.  In particular a tail can first exceed $k$ only if its
core has degree exactly $k$.
\end{lemma}
\begin{proof}
Bottom trimming does not increase degree.  For a retained target, replace the
initial part of an old source path by the horizontal segment from the new
source to its first point in the new bottom row; this does not add a bend.
Conversely, a target below that row is reached with at most one bend, while a
retained target is reached by first going north to the new source and then
following a trimmed path.  Hence a trim lowers degree by at most one.
These claims use Lemma~\ref{pf:source-distance}.

A source-column extension raises degree by at most one, by the half-turned
form of \eqref{pc:tail-degree-bound}.  Shortening that new column removes
paths and a target of source distance zero, so cannot decrease the maximum.
Finally extension and bottom trimming commute up to translation.  These facts
prove both monotonicities and the three possible values.  The last assertion
also follows directly from \eqref{pc:tail-degree-bound}.
\end{proof}

\begin{lemma}\label{pc:width-survival}
If a trim $Q<r$ removes the old rightmost column of a directed convex core,
the trimmed core has degree at most one.  Consequently no extension sequence which
ultimately exceeds a bound $k\ge2$ has a width-losing cumulative trim.
Completion and such retained-width trims commute, and delete precisely the
forest vertices in the removed rows without changing any surviving parent.
\end{lemma}
\begin{proof}
Once the rightmost column has disappeared, the peak shape of the right row
endpoints makes the new bottom row span the entire retained width.  Since
$Q<r$, its left endpoint is still the source column.  Every target is reached
by going east in that bottom row and then north, giving degree at most one.
A later cumulative trim cannot increase degree, and one source extension
raises it by at most one.  Lemma~\ref{pc:history-contraction} therefore
precludes a later degree above two.

The assertions concerning completion and the forest are
Lemma~\ref{pf:bottom-pruning}.
\end{proof}

\begin{lemma}\label{pc:degree-two-tails}
For a rectangular ascending core of degree two, every nested tail
which has a shadow of degree three is ascending; after its first such shadow,
every nested continuation is ascending.
\end{lemma}
\begin{proof}
The full prefix containing that shadow is convex, and its northeast degree
is at least three by \eqref{pc:shadow-maximum}.  It is ascending by
Lemma~\ref{pc:directional-structure}.  Every nested continuation is convex
and preserves the northeast distances between cells of the prefix: an
eastward path cannot visit a later column and return.  The same argument
therefore makes every continuation ascending.
\end{proof}

On an extension sequence ending at its first passage at a bound $k\ge2$, every shadow retains the core
width by Lemma~\ref{pc:width-survival}.  Every earlier tail column meets row
$b_j$ because it contains $I_j$.  Thus the $j$th shadow has width
$\operatorname{width}(R)+j$ and height $b_j+1$.  With $Q_j=H_0-b_j$,
\begin{equation}\label{pc:augmented-weight}
 x^{\operatorname{sp}(P^{[j]})+Q_j}=x^{\operatorname{sp}(R)+j}.
\end{equation}
Each tail column therefore has weight $x$ when the deleted height is retained
in the exponent.  Notice that $Q_j<r$ follows from $b_j\ge a_0$.

\section{The fixed boundary for the cumulative slacks}\label{sec:fixed-boundary}

Fix $D_0$ of degree $k\ge2$ and source height $r$.  The \emph{slack plane}
has integer coordinates $(P,Q)$.  A pair is \emph{admissible} when
$P,Q\ge0$ and $P+Q<r$.  It is \emph{safe} when $\delta(P,Q)\le k$ and
\emph{forbidden} when $\delta(P,Q)>k$.  A \emph{slack path} is a path from $(0,0)$ in
this plane with unit east and north steps.  Its \emph{diagonal level} is
$P+Q$, equal to its number of steps from $(0,0)$.

Use the forest construction of Proposition~\ref{pf:forest-bijection}.  Its row
markers are denoted $r_y$ and its root depths are one.  A source-row marker
of depth at least two lies below the lowest row of the old rightmost column:
unless the rightmost column consists only of its top cell, that latter row
is a depth-one row marker in the completion, and depths are nonincreasing
upwards.  In the exceptional case the rightmost column plainly survives every
trim below any row marker.  Thus every trim below a retained marker of depth at
least two preserves width.

\begin{proposition}\label{pc:fixed-wall}\label{prop:stretched-cut-wall}
Let $D_0$ have degree $k\ge2$ and source height $r$.  Set $H=k+1$ and
\[
 \rho(P)=\operatorname{dep}(r_{r-P-1}),\qquad0\le P<r,
\]
\[
 \begin{aligned}
 a&=\min(\{P:\rho(P)\ge k\}\cup\{r\}),\\
 b&=\min(\{P:\rho(P)=H\}\cup\{r\}).
 \end{aligned}
\]
Then $0\le a\le b\le r$ and $b\ge1$.

Read the cut word $\mu$ backwards for even $k$, forwards for odd $k$.
Replace each letter of $F_e$ for even $k$ (of $F_s$ for odd $k$) by as many
east steps as its tree has depth-$k$ row markers in source-column rows.
Replace each opposite-forest letter by as many north steps as it has depth-$H$
row markers in those rows.  Write $\Gamma$ for the resulting path and call
$\Omega=E^a\Gamma$ the \emph{wall}, starting at $(0,0)$.  Then $\Omega$ ends at $(b,r-b)$ and has length $r$.
For $P<b$, let $\tau(P)$ be the height of its east step at abscissa $P$.
For every admissible $(P,Q)$,
\begin{equation}\label{pc:wall-test}
 \delta(P,Q)=k+1\quad\Longleftrightarrow\quad
 \bigl(P<b\text{ and }Q<\tau(P)\bigr)\ \text{or}\ P\ge b.
\end{equation}
In particular all points with $P=0$ are safe.
\end{proposition}
\begin{proof}
Both forest heights are at most $H$ by
Propositions~\ref{pf:bounce-degree} and~\ref{pf:completion-degree}.
The source column stops below the global top, since otherwise its source
reaches every target with at most one bend.  Its top cell is a column vertex
whose parent is the row marker in that row.  Hence $\rho(0)+1\le H$.
Lemma~\ref{pf:row-order} makes $\rho$ nondecreasing, so the source rows consist of
$a$ shallow rows, $b-a$ depth-$k$ rows and $r-b$ depth-$H$ rows.  The latter
are the entire bottom depth-$H$ row block; no deeper row can lie above the
source top.  This proves the endpoint and length assertions.

If a considered trim loses width, Lemma~\ref{pc:width-survival} gives
$\delta(P,Q)\le2$.  Since trimming lowers degree by at most one, this cannot
happen for $k\ge3$.  For $k=2$ its retained prospective parent must have depth
one, by the observation preceding the proposition, and hence $P<a$ and
$\tau(P)=0$.  Thus \eqref{pc:wall-test} holds in this case.
Otherwise the row markers and their ancestor chains survive the trim by
Lemma~\ref{pf:bottom-pruning}.  By Lemma~\ref{pf:source-extension}, with the
incoming one-step tails in the bounce-path definition, the new source
leaf and its ancestors form the vertical-first chain, and the bottom retained
row marker and its ancestors form the horizontal-first chain.  The latter
has length at most $H$.

When $P<a$, the first chain has length at most $k$, and
Proposition~\ref{pf:bounce-degree} gives associated degree at most $k$.
When $P\ge b$, the new leaf raises its forest to height $H+1$, while the
opposite forest stays at most $H$.  The associated parallelogram is nonflat
of degree $H$, and its directed degree is $H$ by
Proposition~\ref{pf:completion-degree}.

Suppose $a\le P<b$ and put $m=r-b$.  The vertical-first chain has length $H$.
If $Q\ge m$, the horizontal-first chain is shorter, giving degree at most
$k$.  If $Q<m$, the two chains have length $H$: the completion of the extended polyomino is flat of degree $H$.  Its critical tree containing the new leaf is the tree of the selected parent;
its opposite critical tree is the tree of the bottommost retained depth-$H$
row marker.  The pruned rows have depth $H\ge3$, so all roots and their cut
positions survive.

For even $k$, the opposite critical root is the last surviving $F_s$ root
with depth-$H$ row markers.  For odd $k$ it is the first such $F_e$ root.
Lemma~\ref{pf:critical-order} therefore says that the extended polyomino has degree $k+1$ precisely
when a retained opposite depth-$H$ row belongs to a root after the selected
letter in the first case, or before it in the second.  By the row-block order
in Lemma~\ref{pf:row-order}, those rows form a bottom segment of the
terminal depth-$H$ block.  Its size is exactly $\tau(P)$, the number of north
steps preceding the selected east step in $\Gamma$.  Deleting $Q$ bottom rows
leaves such a row exactly when $Q<\tau(P)$, proving the strict cutoff.
Finally $\operatorname{Ext}_0$ duplicates the source column, so
Lemma~\ref{pc:history-contraction} and the elementary bar cases give
$\delta(0,Q)\le k$.
\end{proof}

\begin{corollary}\label{pc:contact-rule}
Read the wall and a slack path at equal diagonal level smaller than $r$.
If the wall abscissa
exceeds the slack abscissa, the slack point is safe.  At a common vertex with
remaining wall length positive, the slack point is safe exactly when the next
wall step is east.
\end{corollary}
\begin{proof}
At the same diagonal level a point strictly to the left of the monotone wall
lies above its horizontal cutoff, hence is safe by \eqref{pc:wall-test}.
At a common vertex $(P,Q)$ followed by an east wall step one has
$Q=\tau(P)$, which is safe.  If the next step is north, the next east step
at abscissa $P$, if any, has height greater than $Q$, giving
$Q<\tau(P)$; if there is no such east step then $P\ge b$.  Both cases are
forbidden.  The positive remaining length ensures admissibility.
\end{proof}

\begin{proposition}\label{pc:canonical-history-weight}
Fix a core and give each added column weight $x$.  A canonical
extension sequence which first hits the forbidden set
\eqref{pc:wall-test} corresponds to a northeast unit path in the slack plane:
replace a column increment $(p_i,q_i)$ by $N^{q_i}E^{p_i}$, and stop the last
east run at its first forbidden point.  This path begins with $N$ and stays safe before its last east step.
For a fixed such path, sum $x^m$ over all extension sequences of $m$ columns
that yield this path and whose final increment stops at its first forbidden
point.  The resulting generating function is
\begin{equation}\label{eq:turn-weighted-history}
 \frac{x^{c+1}}{(1-x)^{s-1}},
\end{equation}
where $s$ is its number of unit steps and $c$ its number of adjacent $EN$
pairs.  If its endpoint is $(P_*,Q_*)$ and $n_*=r-P_*-Q_*$, the final column
may overshoot horizontally by $0,\ldots,n_*-1$; with $u$ marking the final column length, its generating
polynomial is $u+\cdots+u^{n_*}$.
\end{proposition}
\begin{proof}
Increasing $Q$ preserves safety, and increasing $P$ at fixed $Q$
preserves membership in the forbidden set, by Lemma~\ref{pc:slack-monotonicity}.  This proves the unit-path
statement, and $q_0\ge1$ gives the initial north step.
Delete zero-increment columns.  Decomposing the unit word into nonempty
$N^*E^*$ blocks requires a break at every $EN$, while all other gaps are
optional.  With $L$ blocks there may be any number of zero-increment columns
in the $L-1$ intervening gaps, but none before the first canonical column or
after the first forbidden column.  Their combined weight is
$x^L/(1-x)^{L-1}$.  Summing the optional breaks gives
\[
 x\left(\frac{x}{1-x}\right)^c
 \left(1+\frac{x}{1-x}\right)^{s-1-c}
 =\frac{x^{c+1}}{(1-x)^{s-1}}.
\]
The final overshoot shortens the positive source height by any amount less
than $n_*$ without adding another column, giving the stated polynomial.
\end{proof}

\begin{figure}[p]
\normalfont\normalsize
\centering
\begin{minipage}{\linewidth}
\setlength{\parindent}{0pt}
\setlength{\parskip}{5pt}
\begin{tikzpicture}[x=1cm,y=1cm,font=\small]
\path[use as bounding box] (0,0) rectangle (15.6,19.6);

\node[anchor=west,font=\bfseries\normalsize] at (.05,19.18) {A. Source rows and wall steps};
\node[anchor=east] at (15.5,19.18) {$k=3,\quad r=3$};
\node[anchor=west,font=\footnotesize] at (.05,18.57)
  {Source rows of $D_0$, from top to bottom:};
\fill[coreblue!65] (.05,16.78) rectangle (7.1,17.42);
\fill[ghost!65] (.05,15.47) rectangle (7.1,16.78);
\draw[ink!45] (.05,17.49)--(7.1,17.49);
\draw[ink!45] (.05,15.42)--(7.1,15.42);
\node[font=\footnotesize] at (.5,17.95) {$P$};
\node[font=\footnotesize,align=center] at (1.67,17.95) {source row\\$y$};
\node[font=\footnotesize,align=center] at (3.12,17.95) {depth\\$\rho(P)$};
\node[font=\footnotesize] at (4.8,17.95) {forest};
\node[font=\footnotesize,align=center] at (6.28,17.95) {wall\\step};
\foreach \yy/\pp/\rr/\dd/\ff/\ss in
  {17.1/0/2/3/F_s/E,16.45/1/1/4/F_e/N,15.8/2/0/4/F_e/N}{
  \node at (.5,\yy) {$\pp$};
  \node at (1.67,\yy) {$\rr$};
  \node at (3.12,\yy) {$\dd$};
  \node at (4.8,\yy) {$\ff$};
  \node[text=wallblue] at (6.28,\yy) {$\ss$};
}
\node at (3.55,14.95) {$\rho=(3,4,4),\qquad a=0,\quad b=1$};
\node[anchor=west,font=\footnotesize] at (7.8,18.57) {The cut word assigns the steps:};
\node[align=center] at (11.6,17.46) {
  $\displaystyle\mu=\underbrace{S_1}_{F_s}\,
       \underbrace{E_2}_{F_e}\ \longmapsto\
       \underbrace{E}_{\text{one row}}\,
       \underbrace{NN}_{\text{two rows}}$};
\node[draw=wallblue!45,fill=coreblue!45,rounded corners=2pt,inner sep=6pt]
  at (11.6,16.28) {$\displaystyle\Omega=E^a\Gamma=ENN$};
\node[anchor=north west,text width=7.4cm,align=left,font=\footnotesize]
  at (7.8,15.55) {The letter $E_2$ occurs once; its two depth-four
  source rows contribute two north steps.};
\node[anchor=north west,text width=15.15cm,align=left,font=\footnotesize]
  at (.05,14.36) {There is no source row of depth below $k$, so $a=0$.
  The first row of depth $k+1$ has index $P=1$, so $b=1$.};
\draw[ink!25] (0,13.55)--(15.6,13.55);

\node[anchor=west,font=\bfseries\normalsize] at (.05,13.08) {B. Nested intervals};
\node[anchor=west,font=\footnotesize] at (.05,12.5)
  {Original ascending seam and tail; row indices $3,4,5$.};
\foreach \yy/\rr in {9.45/3,10.25/4,11.05/5}
  \node[font=\footnotesize] at (.4,\yy) {$\rr$};
\foreach \xx/\ii/\range in {1.45/0/{[3,5]},3.8/1/{[3,4]},6.15/2/{[4,4]}}{
  \node[align=center] at (\xx,11.92) {$I_{\ii}$\\$\range$};
  \foreach \yy in {9.05,9.85,10.65}
    \draw[cell,fill=ghost!45,dashed] (\xx-.52,\yy) rectangle +(1.04,.8);
}
\foreach \yy in {9.05,9.85,10.65}
  \draw[cell,fill=coreblue] (.93,\yy) rectangle +(1.04,.8);
\foreach \yy in {9.05,9.85}
  \draw[cell,fill=tailorange!55] (3.28,\yy) rectangle +(1.04,.8);
\draw[cell,fill=tailgreen!55] (5.63,9.85) rectangle +(1.04,.8);
\draw[ink!65,-{Stealth[length=2mm]}]
  (2.02,11.45)--(2.57,11.45)--(2.57,10.65)--(3.22,10.65);
\node[font=\footnotesize,anchor=south] at (2.58,11.48) {$q_0=1$};
\draw[ink!65,-{Stealth[length=2mm]}]
  (4.37,9.05)--(4.96,9.05)--(4.96,9.85)--(5.57,9.85);
\node[font=\footnotesize,anchor=south] at (4.96,10.02) {$p_1=1$};
\node[font=\footnotesize] at (1.45,8.6) {seam};
\node[font=\footnotesize,align=center] at (3.8,8.47) {first new\\column};
\node[font=\footnotesize,align=center] at (6.15,8.47) {second new\\column};
\node[align=center,font=\footnotesize] at (11.55,12.03)
  {$I_j=[a_j,b_j]$\\[3pt]$P_j=a_j-a_0,\qquad Q_j=b_0-b_j$};
\node[align=center,font=\footnotesize] at (11.55,10.33) {
  $\renewcommand{\arraystretch}{1.45}
  \begin{array}{c|c|c}
   j&(p_{j-1},q_{j-1})&(P_j,Q_j)\\\hline
   0&\text{--}&(0,0)\\
   1&(0,1)&(0,1)\\
   2&(1,0)&(1,1)
  \end{array}$};
\node[draw=pathred!45,fill=pathred!4,rounded corners=2pt,inner sep=6pt]
  at (11.55,8.55) {$\displaystyle\gamma=N^{q_0}E^{p_0}N^{q_1}E^{p_1}=NE$};
\node[anchor=west,font=\footnotesize] at (.05,7.8)
  {Here each new column gives one unit step.};
\draw[ink!25] (0,7.35)--(15.6,7.35);

\node[anchor=west,font=\bfseries\normalsize] at (.05,6.88) {C. The fixed boundary};
\node[anchor=east] at (15.5,6.88) {$P,Q\ge0,\qquad P+Q<3$};
\begin{scope}[shift={(1.12,1.03)},x=1.55cm,y=1.55cm]
  \fill[ghost!80] (0,3)--(3,3)--(3,0)--cycle;
  \fill[safeolive!3] (0,0)--(3,0)--(0,3)--cycle;
  \draw[step=1,ink!12,line width=.25pt] (0,0) grid (3,3);
  \draw[axis] (0,0)--(3.3,0) node[anchor=west] {$P$};
  \draw[axis] (0,0)--(0,3.3) node[anchor=south] {$Q$};
  \draw[ink!55,densely dashed] (0,3)--(3,0);
  \node[rotate=-45,fill=ghost!80,inner sep=2pt,font=\footnotesize]
    at (2.08,1.25) {$P+Q=3$};
  \draw[wall] (0,0)--(1,0)--(1,2);
  \node[text=wallblue] at (.52,.18) {$E$};
  \node[text=wallblue] at (1.18,.43) {$N$};
  \node[text=wallblue] at (1.18,1.52) {$N$};
  \draw[slack,-{Stealth[length=2mm]}] (0,0)--(0,1)--(1,1);
  \node[text=pathred] at (-.18,.5) {$N$};
  \node[text=pathred] at (.46,1.17) {$E$};
  \foreach \xx/\yy in {0/0,0/1,0/2}
    \fill[safeolive] (\xx,\yy) circle[radius=2.35pt];
  \foreach \xx/\yy in {1/0,1/1,2/0}{
    \draw[pathred,line width=1.25pt]
      (\xx-.05,\yy-.05)--(\xx+.05,\yy+.05)
      (\xx-.05,\yy+.05)--(\xx+.05,\yy-.05);
  }
  \draw[pathred,line width=.7pt] (1,1) circle[radius=4.6pt];
  \foreach \xx/\yy in {0/3,1/2,2/1,3/0}
    \draw[ink!55,fill=white,line width=.6pt] (\xx,\yy) circle[radius=2.7pt];
  \node[anchor=north east,font=\footnotesize] at (-.01,-.06) {$(0,0)$};
  \node[anchor=east,font=\footnotesize] at (-.08,1) {$(0,1)$};
  \node[anchor=east,font=\footnotesize] at (-.08,2) {$(0,2)$};
  \node[anchor=north,font=\footnotesize] at (1,-.11) {$(1,0)$};
  \node[anchor=north,font=\footnotesize] at (2,-.11) {$(2,0)$};
  \node[anchor=west,font=\footnotesize,text=pathred] at (1.12,1.03) {$(1,1)$};
  \node[anchor=west,font=\footnotesize] at (1.1,2.15) {$(1,2)$: end of $\Omega$};
\end{scope}
\draw[wall] (8.1,6.07)--(8.8,6.07);
\node[anchor=west] at (8.94,6.07) {$\Omega=ENN$};
\draw[slack] (8.1,5.54)--(8.8,5.54);
\node[anchor=west] at (8.94,5.54) {$\gamma=NE$};
\node[anchor=north west,text width=7.15cm,align=left,font=\footnotesize]
  at (8.1,4.99) {$P$: lower endpoint moves up.\\$Q$: upper endpoint moves down.};
\node[anchor=north west,text width=7.15cm,align=left,font=\footnotesize]
  at (8.1,3.99) {The canonical first step is $N$, from $(0,0)$ to $(0,1)$.};
\node[anchor=north west,text width=7.15cm,align=left,font=\footnotesize]
  at (8.1,2.98) {At $(1,1)$ the next wall step is $N$: this is the first
  forbidden point of $\gamma$.};
\node[anchor=north west,text width=7.15cm,align=left,font=\footnotesize]
  at (8.1,1.7) {The terminal diagonal $P+Q=3$ is excluded.};
\fill[safeolive] (.5,.25) circle[radius=2.35pt];
\node[anchor=west,font=\footnotesize] at (.72,.25) {safe};
\draw[pathred,line width=1.25pt] (4,.2)--(4.1,.3) (4,.3)--(4.1,.2);
\node[anchor=west,font=\footnotesize] at (4.25,.25) {forbidden};
\draw[ink!55,fill=white,line width=.6pt] (8.35,.25) circle[radius=2.7pt];
\node[anchor=west,font=\footnotesize] at (8.58,.25) {inadmissible};
\end{tikzpicture}
\par
\caption{\textbf{Nested intervals and the fixed boundary.} The directed core from Figures~\ref{fig:core-shadow}--\ref{fig:completion-cut} has rows $[0,1],[0,1],[0,2],[1,3],[2,4],[3,3]$
from bottom to top.  Its three source rows have depths $3,4,4$ from top to bottom.
For $k=3$, the cut word $S_1E_2$ therefore gives the boundary $ENN$.
The original seam interval $[3,5]$ is shortened first from above and then from below,
giving the canonical path $NE$.  In general, one added column with increments $(p,q)$
is represented by $N^qE^p$, which may have several unit steps; here the two columns
each give one step. The TikZ code for this figure was prepared with assistance from ChatGPT (OpenAI) and reviewed by the author.}
\label{fig:fixed-boundary}
\end{minipage}
\end{figure}
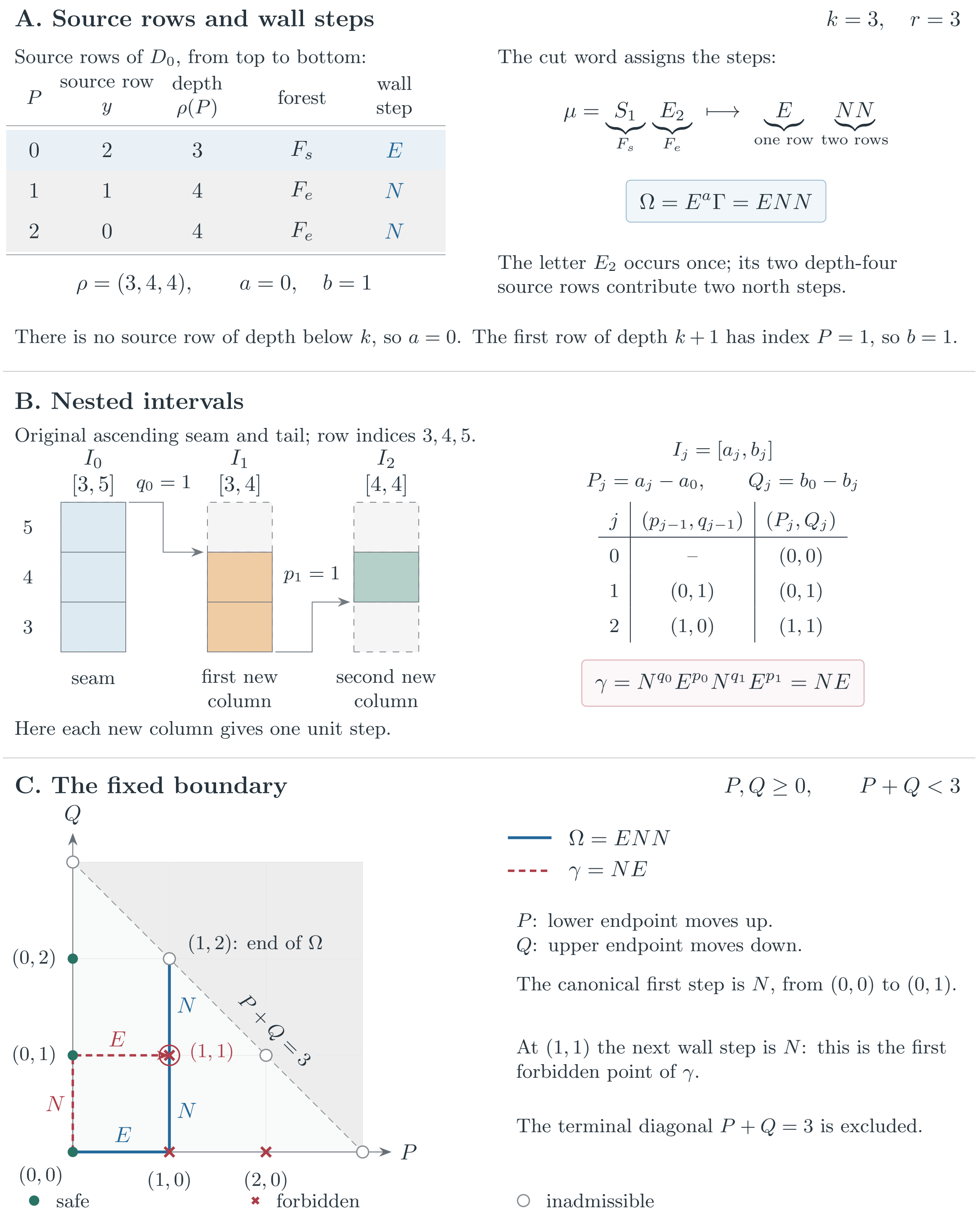

Figure~\ref{fig:fixed-boundary} follows the example of Figures~\ref{fig:core-shadow} and~\ref{fig:completion-cut}: the source-row depths give the wall $ENN$, and the two nested columns give the slack path $NE$, whose endpoint is its first forbidden point.

\section{Separation at the first forbidden shadow}\label{sec:separation}

Let $\mathcal T_n$ count arbitrary nested column sequences appended to an interval of
length $n$, including the empty sequence, each new column carrying weight $x$.
Let $\widehat{\mathcal T}_n$ impose, for a nonempty sequence, avoidance of the top cell
by its first column.  Their series will be evaluated in Section~\ref{pe:enumeration}.

Let $B_{k,n}(x)$ be the generating function of canonical prefixes
$R,I_1,\ldots,I_m$ with $m\ge1$, $\deg R=k$, last column length $n$,
$\deg P^{[i]}\le k$ for $0\le i<m$, and $\deg P^{[m]}=k+1$.
Each prefix has weight $x^{\operatorname{sp}(R)+m}$, its actual
semi-perimeter weight by Lemma~\ref{pc:core-tail}.  We abbreviate this
series to $B_{k,n}$.  Let $n(R)$ be the height of the seam column of $R$.
For every $k\ge2$,
\begin{equation}\label{thm:first-passage-factorisation}
 A-A_k=
 \sum_{\substack{R\ \text{rectangular ascending}\\\deg R>k}}
 x^{\operatorname{sp}(R)}\widehat{\mathcal T}_{n(R)}
 +\sum_{n\ge1}B_{k,n}\mathcal T_n,
\end{equation}

To prove the identity, a core of degree greater than $k$ has degree at least three, so all
nested tails are allowed by Lemma~\ref{pc:free-tails}.  Otherwise
\eqref{pc:shadow-maximum} locates a unique first forbidden shadow.  Trimming lowers
degree and one column raises it by at most one, so its degree is $k+1$.
Equation~\eqref{pc:tail-degree-bound} forces the core degree to be $k$.
For $k\ge3$ all nested continuations are ascending.  For $k=2$ exactly the
same conclusion, including validity of every first-passage prefix,
is Lemma~\ref{pc:degree-two-tails}.  Thus both products count precisely the
indicated disjoint classes.  This also proves that every extension sequence used in
Proposition~\ref{pc:canonical-history-weight} satisfies the geometric and
weight conditions needed for the factorization.

\section{Enumeration of the excluded ascending class}
\label{pe:enumeration}

Throughout this section $x$ marks semi-perimeter.  We use the forest and
source-row descriptions established in Section~\ref{pf:forest-section}, and the
canonical core decomposition of Section~\ref{pc:core-section} and the
first-passage identity \eqref{thm:first-passage-factorisation}.
We take $k\ge2$.  All series below are formal power series; every infinite
sum is coefficientwise finite because its objects have bounded size at each
power of $x$.

\subsection{Trees with a distinguished source row}
The bounded-height tree series is the usual Fibonacci-polynomial expression;
see~\cite[Section~2 and Proposition~1]{BRSp}.  We record its recurrence to fix the height
index before introducing the source-row markers.

Put
\[
 \Phi_0=0,\qquad \Phi_1=1,\qquad
 \Phi_{j+1}=\Phi_j-x\Phi_{j-1},\qquad
 T_{\le j}=\frac{x\Phi_j}{\Phi_{j+1}}.
\]
For a tree $T$, let $|T|$ be its number of vertices, $\operatorname{ht}(T)$
its height, and $n_j(T)$ its number of vertices at depth $j$.  Precisely,
\[
 T_{\le j}=\sum_{\operatorname{ht}(T)\le j}x^{|T|},\qquad
 T_j(u)=\sum_{\operatorname{ht}(T)\le j}x^{|T|}u^{n_j(T)}.
\]
Thus both series restrict the height to at most $j$, and $T_{\le0}=0$.
A root contributes $x$;
its ordered sequence of children gives
$T_{\le j}=x/(1-T_{\le j-1})$.  The displayed polynomial recurrence solves
this identity with its empty initial class.  Iteration of $f(z)=x/(1-z)$ gives
\begin{equation}\label{pe:tree-marking}
 T_j(u)=\frac{x(\Phi_{j-1}-ux\Phi_{j-2})}
                   {\Phi_j-ux\Phi_{j-1}}
 =T_{\le j-1}+\frac{U_j u}{1-T_{\le j-1}u},
 \qquad U_j=\frac{x^j}{\Phi_j^2},
\end{equation}
for $j\ge2$, while $T_1(u)=xu$.  The second equality follows from
$\Phi_{j-1}^2-\Phi_{j-2}\Phi_j=x^{j-2}$, which itself follows by induction
from the recurrence.  Thus the trees of height at most $j$ with exactly $m\ge1$ vertices
at depth $j$ have generating function $U_jT_{\le j-1}^{m-1}$.

The source parent is the last non-leaf in its row-marker level in the
source tree; the row markers after it are the same-tree source rows below
it.  This is the source extremality statement of
Lemma~\ref{pf:source-extremal}.  For $j\ge2$ define
\begin{align}
 \Sigma_j(u)&=
 \frac{x^j}{\Phi_{j+1}(\Phi_{j-1}-ux\Phi_{j-2})},
 \label{pe:sigma}\\
 \Sigma_j^+(u,v)&=
 \frac{x^j}{(\Phi_{j-1}-ux\Phi_{j-2})
                 (\Phi_{j+1}-vx\Phi_j)}.
 \label{pe:sigma-plus}
\end{align}
The first counts source trees whose source top has depth $j$ and
whose tree height is $j$;
$u$ counts depth-$(j-1)$ source rows strictly below the parent.
The second permits tree height at most $j+1$, including height $j$,
still with source-top depth $j$;
$v$ counts all depth-$(j+1)$ row markers.  We abbreviate
$\Sigma_j^+(u,u)$ to $\Sigma_j^+(u)$.

Here is a direct derivation.  At row depth $j-1$, the generating function for a subtree of height at most two is
$J=x/(1-x)$.  Heights within this subtree are measured from its own root.  The distinguished last non-leaf has weight $J-x$;
subsequent row markers are leaves of weight $ux$.  Consider
\[
 M=\begin{pmatrix}J&J-x\\0&ux\end{pmatrix}
\]
and iterate $f(M)=x(I-M)^{-1}$ another $j-2$ times; $I$ denotes the
identity matrix.  We write $f^{\circ m}$ for $m$-fold composition, with
$f^{\circ0}(M)=M$.  The diagonal entries describe the subtrees before and
after the distinguished parent at that level.  Each contribution to the
upper-right entry of a product of upper triangular matrices contains
exactly one upper-right factor, which selects that parent.  The sum
$(I-M)^{-1}$ permits any number of ordered children.
Thus the upper-right entry of $f^{\circ(j-2)}(M)$ counts the desired source tree.  For the second refinement use instead
\[
 J=J(v)=\frac{x(1-vx)}{1-x-vx}.
\]
This counts subtrees of height at most three, with $v$ recording
the number of vertices at depth three relative to the subtree root.
For a triangular matrix the upper-right entry of $f^{\circ m}(M)$ is
\[
 (J-x)\,
 \frac{x^m}{(\Phi_{m+1}-J\Phi_m)(\Phi_{m+1}-ux\Phi_m)}.
\]
Taking $m=j-2$ gives \eqref{pe:sigma}.  In the height-three case use
\[
 J(v)-x=\frac{x^2}{1-x-vx},\qquad
 (1-x-vx)(\Phi_{j-1}-J(v)\Phi_{j-2})
 =\Phi_{j+1}-vx\Phi_j
\]
to obtain \eqref{pe:sigma-plus}.  The last child of the distinguished parent
is uniquely the source top, so there is no further choice of a
distinguished vertex.  All rows
of depth $j+1$ lie below it and belong to the source rows, including any
leaves below the source-top vertex itself.  The source parent's own row is
counted separately.

\subsection{Cores of a prescribed degree}

Let $\mathcal C_q(u)$ count directed convex polyominoes of degree exactly
$q$, with $x$ marking semi-perimeter and $u$ marking source height.  These
are the half-turned rectangular ascending cores.  For $q\ge2$, let
$N_q^{\mathrm{hi}}(u)$ and $N_q^{\mathrm{lo}}(u)$ count those with nonflat
completion of degree $q$, according as the source top has depth $q+1$ or
$q$.  The cut contributes the factor $x^2$.  We have
\begin{align}
 N_q^{\mathrm{hi}}(u)
 &=\frac{ux^2\Sigma_{q+1}(u)}
 {(1-T_q(u)-T_{\le q})(1-T_{\le q+1}-T_{\le q})},
 \label{pe:Nhi}\\
 N_q^{\mathrm{lo}}(u)
 &=\frac{ux^2\Sigma_q^+(u)}
 {(1-T_{q-1}(u)-T_q(u))(1-T_{q+1}(u)-T_q(u))}
 \nonumber\\[-1mm]
 &\quad-
 \frac{ux^2\Sigma_q(u)}
 {(1-T_{q-1}(u)-T_q(u))(1-T_{\le q}-T_q(u))}.
 \label{pe:Nlo}
\end{align}
For flat completions of flat height $h\ge2$, let $F_h^{\mathrm{all}}(u)$
count all cuts and $F_h^{\mathrm{av}}(u)$ those avoiding the critical
rectangle, with the same variables $x,u$.  Use the abbreviations
\[
 A=T_{\le h-1},\quad A_u=T_{h-1}(u),\quad
 B=T_{\le h},\quad B_u=T_h(u)
\]
and define
\begin{align}
 F_h^{\mathrm{all}}(u)
 &=ux^2\Sigma_h(u)
 \left\{\frac1{(1-A_u-B_u)(1-B-B_u)}
       -\frac1{(1-A_u-A)(1-B-A)}\right\},
 \label{pe:Fall}\\
 F_h^{\mathrm{av}}(u)
 &=\frac{ux^2\Sigma_h(u)(B_u-A)}
 {(1-A_u-B_u)(1-A_u-A)(1-B-A)}.
 \label{pe:Fav}
\end{align}
For $q\ge3$, the complete degree class is
\begin{equation}\label{pe:source-layer}
 \mathcal C_q(u)=N_q^{\mathrm{hi}}(u)+N_q^{\mathrm{lo}}(u)
 +F_q^{\mathrm{all}}(u)-F_q^{\mathrm{av}}(u)
 +F_{q+1}^{\mathrm{av}}(u).
\end{equation}

To prove these formulas, first suppose the forest heights are unequal, with
maximum $q+1$.  By Proposition~\ref{pf:bounce-degree}, the source top has depth
$q+1$ or $q$.  In the former case, trees of the source forest preceding the source
tree have height at most $q$, and those following it have height at most
$q+1$; all trees of the other forest have height at most $q$.  Source-row order gives the two factors in
\eqref{pe:Nhi}, with $\Sigma_{q+1}$ for the distinguished tree.  In the
latter case the source forest has height $q+1$, its source top has depth
$q$, and the source parent has depth $q-1$.  Trees of the source forest preceding the source tree have height at
most $q-1$, and those following it have height at most $q+1$.
The opposite forest has bound $q$.  This gives the first term in
\eqref{pe:Nlo}; its second term is precisely the subfamily in which the
source forest never reaches $q+1$.  The factor $u$ in both formulas counts
the source parent's row.  For one parity the source root is first among the
roots reaching its depth; for the other it is last.  Reversal of the cut word and exchange of the forest names give the same products.

The two nonflat formulas include $q=2$.  The source-top depths are then
three and two; in the latter case its parent has depth one and is a root.
The same order argument applies with $T_1(u)=xu$.  The source-tree
calculation uses zero iterations of $f$ and gives
\[
 \Sigma_2(u)=\frac{x^2}{1-x},\qquad
 \Sigma_2^+(u,v)=\frac{x^2}{1-x-vx}.
\]
The subtraction in $N_2^{\mathrm{lo}}$ removes exactly the forests whose
source forest has height at most two.

For the flat case take odd $h$, so the source tree is the first
height-$h$ tree of $F_e$.  Its weight is $\Sigma_h(u)$.  The preceding $F_e$ roots carry
$A_u$, the following ones $B$, and every $F_s$ tree carries $B_u$ because
all depth-$h$ row markers lie in the source rows.  Requiring at least one
height-$h$ tree of $F_s$ gives \eqref{pe:Fall} by subtraction.  By
Lemma~\ref{pf:critical-order}, avoidance requires the root of the last height-$h$ tree of $F_s$ to occur before the source root in the cut word.  The prefix
therefore contains a height-$h$ tree of $F_s$ and the suffix contains none, giving
\[
 ux^2\Sigma_h(u)
 \left(\frac1{1-A_u-B_u}-\frac1{1-A_u-A}\right)\frac1{1-B-A},
\]
which is \eqref{pe:Fav}.  Reversal gives even $h$.  Finally
Proposition~\ref{pf:completion-degree} partitions degree $q$ into the two
nonflat classes, the flat height-$q$ meeting class and the flat
height-$(q+1)$ avoiding class.  This proves \eqref{pe:source-layer}.

\subsection{Canonical tails and the telescoping sum}

Recall the tail series $\mathcal T_n$ and $\widehat{\mathcal T}_n$ from the preceding
decomposition, including the empty tail.
A subinterval of length $j$ has $n-j+1$ placements, of which $n-j$ avoid
that top cell.  Therefore
\[
 \mathcal T_n=1+x\sum_{j=1}^n(n-j+1)\mathcal T_j,
 \qquad
 \widehat{\mathcal T}_n=1+x\sum_{j=1}^{n-1}(n-j)\mathcal T_j.
\]
Multiplication by $z^n$ and summation give
\[
 \sum_{n\ge1}\mathcal T_nz^n=\frac{z(1-z)}{(1-z)^2-x},\qquad
 \sum_{n\ge1}\widehat{\mathcal T}_nz^n=\frac{z(1-z-x)}{(1-z)^2-x}.
\]
Equivalently, with $\rho_\pm=(1\mp\sqrt x)^{-1}$,
\begin{equation}\label{pe:tail-evaluations}
 \mathcal T_n=\tfrac12(\rho_+^n+\rho_-^n),\qquad
 \widehat{\mathcal T}_n=\tfrac12(\rho_+^{n-1}+\rho_-^{n-1}).
\end{equation}
For a series $C(u)=\sum_{n\ge1}c_n(x)u^n$ whose coefficient of each
power of $x$ is a polynomial in $u$, define its \emph{canonical tail sum}
by $\widehat{\mathcal C}[C]=\sum_{n\ge1}c_n(x)\widehat{\mathcal T}_n$.  Then
\begin{equation}\label{pe:canonical-functional}
 \widehat{\mathcal C}[C]
 =\tfrac12\{(1-\sqrt x)C(\rho_+)+(1+\sqrt x)C(\rho_-)\}.
\end{equation}
The two terms are conjugate and their sum belongs to $\mathbb Q[[x]]$.

Let $t$ be the unique series in $x\mathbb Q[[x]]$ satisfying
$t=x(1+t)^2$.  Then
\begin{equation}\label{pe:Catalan}
 x=\frac{t}{(1+t)^2},\qquad t=x(1+t)^2,
 \qquad \Phi_j=\frac{1-t^j}{(1-t)(1+t)^{j-1}}.
\end{equation}
Define
\begin{align}
 M(t,z)&=t^4+t^2+2t
 +(1-t+2t^2+2t^3+t^4-t^5)z
 \nonumber\\
 &\quad-(t+t^3+2t^4)z^2-4t^3z^3,
 \label{pe:M}\\
 R(t,z)&=\frac{t^2zM(t,z)}
 {(1-t)^3(1+t)^5(1-z^2)(1-tz^2)}.
 \label{pe:R}
\end{align}
For $q\ge3$, the series of cores of degree $q$ satisfy
\begin{equation}\label{pe:telescope}
 \widehat{\mathcal C}[\mathcal C_q]
 =R(t,t^q)-R(t,t^{q+1}).
\end{equation}
Here is an explicit algebraic proof, uniform in $q$.  Put $z=t^q$ and
$L=(1+t)^q$.  For every index shift $c$ in \eqref{pe:source-layer}, substitute
\[
 \Phi_{q+c}=\frac{1-zt^c}{(1-t)L(1+t)^{c-1}},\qquad
 x^{q+c}=\frac{zt^c}{L^2(1+t)^{2c}}.
\]
The factors $L$ cancel from each summand.  In
\eqref{pe:canonical-functional} put $t=y^2$ and
$\sqrt x=y/(1+t)$; the arguments become $(1+t)/(1+t\mp y)$.
Substituting the five explicit summands of \eqref{pe:source-layer}, pairing
the signs of $y$, and putting the fractions over a common denominator gives
\begin{equation}\label{pe:telescope-certificate}
 \widehat{\mathcal C}[\mathcal C_q]
 =\frac{t^2zP(t,z)}
 {(1-t)^2(1+t)^4(1-z^2)(1-tz^2)(1-t^2z^2)(1-t^3z^2)},
\end{equation}
where the numerator is explicitly
\begin{align*}
 P={}&t^8z^5+2t^7z^6-t^7z^5-t^7z^4-t^6z^6+2t^6z^5+2t^6z^4\\
 &+t^5z^6+2t^5z^5-t^5z^4-4t^5z^3-t^5z^2-t^5z\\
 &+t^4z^5-t^4z^4-t^4z^2+t^4z-t^3z^5-t^3z^4-4t^3z^3\\
 &-t^3z^2+2t^3z+t^3+2t^2z^2+2t^2z-t^2-tz^2-tz+2t+z.
\end{align*}
For clarity, the numerator identity needed for this denominator clearing and
for the final subtraction is
\begin{equation}\label{pe:P-identity}
 (1-t^2)P(t,z)=M(t,z)(1-t^2z^2)(1-t^3z^2)
 -tM(t,tz)(1-z^2)(1-tz^2).
\end{equation}
Expanding the right side gives the displayed polynomial multiplied by
$1-t^2$.  Putting $R(t,z)-R(t,tz)$ over the same denominator and cancelling
$(1-t)(1+t)$ now gives \eqref{pe:telescope-certificate}, proving
\eqref{pe:telescope}.  This computation is a rational identity in the
independent variables $t,z$, so the single calculation covers every $q\ge3$.

The canonical core contribution to the excluded degree-$k$ class is consequently
\begin{equation}\label{pe:bad-core}
 B_k^{\mathrm{core}}
 =\sum_{q\ge k+1}\widehat{\mathcal C}[\mathcal C_q]
 =R(t,t^{k+1}).
\end{equation}
The last equality follows by telescoping; $R(t,t^q)$ tends coefficientwise
to zero as $q\to\infty$.

\subsection{The two initial wall families}

The fixed-wall theorem (Proposition~\ref{pc:fixed-wall}) shows that only
two parts of the degree-$k$ source
class can have a later first passage: $N_k^{\mathrm{lo}}$ and
$F_{k+1}^{\mathrm{av}}$.  Indeed the other non-flat part has source-parent
depth $k$ but no row of depth $k+1$; a flat height-$k$ core likewise has no
such row.  Its completed wall is horizontal, so cannot be crossed.  In the
two remaining parts the source-row profiles, from top down, are respectively
\[
 (k-1)^A k^L(k+1)^M\quad(A\ge1),
 \qquad k^L(k+1)^M\quad(L\ge1),
\]
Here $d^m$ means $m$ consecutive source rows of depth $d$; $A,L,M$
are local nonnegative row counts and $M\ge1$.  These descriptions follow
from source extremality and row-depth order.

Set $H=k+1$ and
\[
 \varrho=T_{\le k-2},\quad r=T_{\le k-1},\quad
 s=T_{\le k},\quad B=T_{\le H},\quad
 c=U_k,\quad d=U_H,
\]
\[
 S_N=\frac{x^k}{\Phi_{k-1}\Phi_{k+1}},\qquad
 S_F=\frac{x^{k+1}}{\Phi_k\Phi_{k+2}},\qquad
 \mathcal G(u,v)=\frac1{1-T_k(u)-T_H(v)}.
\]
The decomposition below distinguishes an initial segment of east steps.
It need not be the whole initial east run, since the remaining wall may
also start east.  The variable $w$ records this segment's length; $u$ and
$v$ record east and north steps in the remaining wall.  The functions
$P_N(w),P_F(w)$ count the distinguished segment and its associated trees,
including the factor $x^2$ for the cut.  Here $N$ and $F$ denote the nonflat
and flat families.  Explicitly,
\begin{align}
 P_N(w)&=\frac{wx^2S_N}
 {(1-\varrho w)(1-T_{k-1}(w)-T_k(w))},
 \label{pe:PN}\\
 P_F(w)&=\frac{wx^2S_F}
 {(1-rw)(1-T_k(w)-s)(1-B-s)}.
 \label{pe:PF}
\end{align}
The ordered wall series are
\begin{equation}\label{pe:ordered-walls}
 W_N(w;u,v)=P_N(w)\frac1{1-sv}\mathcal G(u,v),\qquad
 W_F(w;u,v)=P_F(w)\frac{dv}{1-sv}\mathcal G(u,v).
\end{equation}
In the first formula the source forest is permitted to have height at most
$H$.  Its exact $N_k^{\mathrm{lo}}$ part is $W_N(w;u,v)-W_N(w;u,0)$.
The removed walls are horizontal and contribute nothing to first passage.
The flat formula is exact without subtraction.

We prove the order of these factors, since their ordinary products alone
would not suffice.  Take odd $k$, so $F_s$ supplies horizontal depth-$k$
rows and $F_e$ vertical depth-$H$ rows.  In the nonflat family the source
tree is the first $F_e$ tree reaching depth $k$.  Before it, the $F_e$
trees are bounded by $k-1$ and contribute their shallow source rows; the
$F_s$ trees contribute horizontal depth-$k$ rows.  There is no vertical
unit in this prefix.  The source tree contributes its remaining shallow
rows, its parent, and then its depth-$H$ vertical block.  Its weight is
\[
 w\Sigma_k^+(w,v)=\frac{wS_N}{(1-\varrho w)(1-sv)}.
\]
All depth-$k$ and depth-$H$ rows after the letter corresponding to the source tree belong to the source rows, and the remaining word is unrestricted, giving $\mathcal G(u,v)$.  The
shallow offset and the pre-source depth-$k$ rows are all horizontal, so their
union has precisely the weight $P_N(w)$.

For the flat family, the source tree is the critical $F_s$ tree $u_s$.
Its cut word has the word decomposition
\[
 \{e_1,s_2\}^*u_s\{e_1,s_1\}^*u_e\{e_2,s_1\}^*.
\]
The first region has no boundary contribution: its $e_1$ trees have no depth-$H$
rows, and the depth-$k$ rows of its $s_2$ roots lie above the source rows.
It contributes $(1-s-B)^{-1}$.  The source tree contributes
$w\Sigma_H(w)=wS_F/(1-rw)$; the middle region contains only horizontal
wall steps and has weight $(1-s-T_k(w))^{-1}$.  The critical $u_e$ tree
then contributes the nonempty north-step sequence
$T_H(v)-s=dv/(1-sv)$, followed by the unrestricted remaining word $\mathcal G(u,v)$.
This is exactly the order asserted in \eqref{pe:ordered-walls}.  Reversal
and forest exchange prove even $k$.  This construction preserves the required order of the east and north
contributions to the wall.

\subsection{Summing over the wall and the slack path}

Let
\[
 Z_0=\frac1{1-r-s},\qquad a=Z_0c,\qquad b=Z_0d.
\]
A tree with no source row at the relevant depth contributes no wall
step.  Such trees have total weight $r+s$; summing arbitrary intervening
sequences gives $Z_0$ initially and the factors $a,b$ at the first steps
of the remaining tree contributions.  By \eqref{pe:tree-marking}, a horizontal run of length
$m\ge1$ therefore has weight $a(r+a)^{m-1}$, and a vertical run has weight
$b(s+b)^{m-1}$.  Since
\[
 s=\frac{x}{1-r},\quad c=(s-r)(1-r),\quad d=s(s-r),
\]
we obtain
\begin{equation}\label{pe:wall-parameters}
 r+a=s+b=\lambda,\qquad
 \lambda=\frac{x\Phi_{2k}}{\Phi_{2k+1}},\qquad
 ab=x-\lambda+\lambda^2.
\end{equation}
The first identity follows by substitution in $a,b$; the other two follow
by the recurrence for $\Phi_j$, or directly by \eqref{pe:Catalan}.
Index matrix rows and columns by $E,N$, in that order.  The entry in row
$i$, column $j$ is the factor for adjoining direction $j$ after direction
$i$.  Sequence weights are products of these factors.  We write $I$ for
the identity matrix of the required size and $\operatorname{diag}(u,v)$
for the diagonal matrix with entries $u,v$.  The weights of successive
wall directions are
\[
 W=\begin{pmatrix}\lambda&b\\a&\lambda\end{pmatrix},
 \qquad
 C_N(u,v)=(0,1)(I-W\operatorname{diag}(u,v))^{-1}\binom11.
\]
Here $C_N(u,v)$ counts the remaining wall, including the empty sequence,
with $u,v$ recording its east and north steps.  Its initial row vector
$(0,1)$ gives the factors for a preceding north direction; it inserts no
extra north step.  Summing alternating runs gives
\begin{equation}\label{pe:north-suffix}
 \frac1{1-sv}\mathcal G(u,v)=Z_0C_N(u,v),\qquad
 \frac{dv}{1-sv}\mathcal G(u,v)=bvC_N(u,v).
\end{equation}
Thus the nonflat remainder uses the initial row vector $(0,1)$;
the flat remainder has one prescribed first north step of weight $b$.

The canonical slack-path formula gives a first north step of weight $x$;
each subsequent step has factor $q=(1-x)^{-1}$ and each $E$-to-$N$ turn
has one additional factor $x$.  Its direction weights are
\[
 V=q\begin{pmatrix}1&x\\1&1\end{pmatrix}.
\]
Compare vertices of the wall and slack path at equal diagonal level.
Their \emph{gap} is wall abscissa minus slack-path abscissa.  For paired
directions $ij$, the first letter is the wall direction and the second
the slack direction.  Use the index order $EE,EN,NE,NN$.
The Kronecker product $W\otimes V$ has entry $W_{ii'}V_{jj'}$ from $ij$
to $i'j'$, so the joint direction factors are
\[
 W\otimes V=q\begin{pmatrix}
 \lambda&\lambda x&b&bx\\
 \lambda&\lambda&b&b\\
 a&ax&\lambda&\lambda x\\
 a&a&\lambda&\lambda
 \end{pmatrix}.
\]
Columns $EN$ and $NE$ respectively increase and decrease the gap;
columns $EE,NN$ preserve it.  Let $M_+$ retain only column $EN$,
$M_-$ only column $NE$, and $M_0$ columns $EE,NN$, setting the other
columns to zero.  Let $e_{ij}$ be the coordinate column vector at index
$ij$; a superscript $T$ denotes transpose.  Let $\mathbf u,\mathbf v$ be columns $EN,NE$ of $W\otimes V$.
Then $M_+=\mathbf u e_{EN}^T$ and $M_-=\mathbf v e_{NE}^T$.

The component $d_{ij}$ of $\mathbf d$ counts paired continuations with
preceding directions $ij$ whose gap first drops from $g$ to $g-1$ on
their last step, staying at least $g$ before it.  Here $g\ge1$, and
translation of the gap shows that the series is independent of $g$.
First-step
and first-return decomposition gives the vector equation
\begin{equation}\label{pe:descent-equation}
 \mathbf d=\mathbf v+M_0\mathbf d+
 \mathbf u(e_{EN}^T\mathbf d)(e_{NE}^T\mathbf d).
\end{equation}
All entries of the joint step matrix have positive $x$-order, so this equation
has a unique formal solution.  It is
\begin{equation}\label{pe:descent}
 \mathbf d=\begin{pmatrix}\chi/(1+t)\\\chi\\t\\t(1+t)\end{pmatrix},
 \qquad \chi=\frac{t-(1+t)\lambda}{a}.
\end{equation}
Although the quotient contains $a$, its numerator has sufficiently high
order, as is also explicit below.  To prove the displayed solution, set $h=1+t$.
The relations
\[
 a\chi=t-h\lambda,\qquad
 bh^2=(1-h\lambda)\chi
\]
follow from \eqref{pe:wall-parameters} and $t=xh^2$.  Substitution of the
candidate vector into the right side of \eqref{pe:descent-equation} gives
\[
 (bh+\lambda\chi,\ bh^2+h\lambda\chi,
   h\lambda+a\chi,\ h^2\lambda+ha\chi)^T,
\]
which is the candidate itself.  This proves \eqref{pe:descent}.
The one-level descent matrix is $D=\mathbf d e_{NE}^T$; hence
\begin{equation}\label{pe:descent-powers}
 D^g=t^{g-1}D\quad(g\ge1).
\end{equation}

\subsection{Contact with the wall and completion of the suffix}

\begin{figure}[p]
\normalfont\normalsize
\centering
\begin{minipage}{\linewidth}
\setlength{\parindent}{0pt}
\setlength{\parskip}{5pt}
The wall $\Omega$ and the slack path $\gamma$ use unit east ($E$) and north ($N$) steps.
Compare them at equal diagonal levels $\ell=P+Q$; the gap is the wall abscissa minus
the slack abscissa.  In a paired direction, the wall step is written first.

\begin{minipage}[t]{\linewidth}
\textbf{(a) The two paths}\par
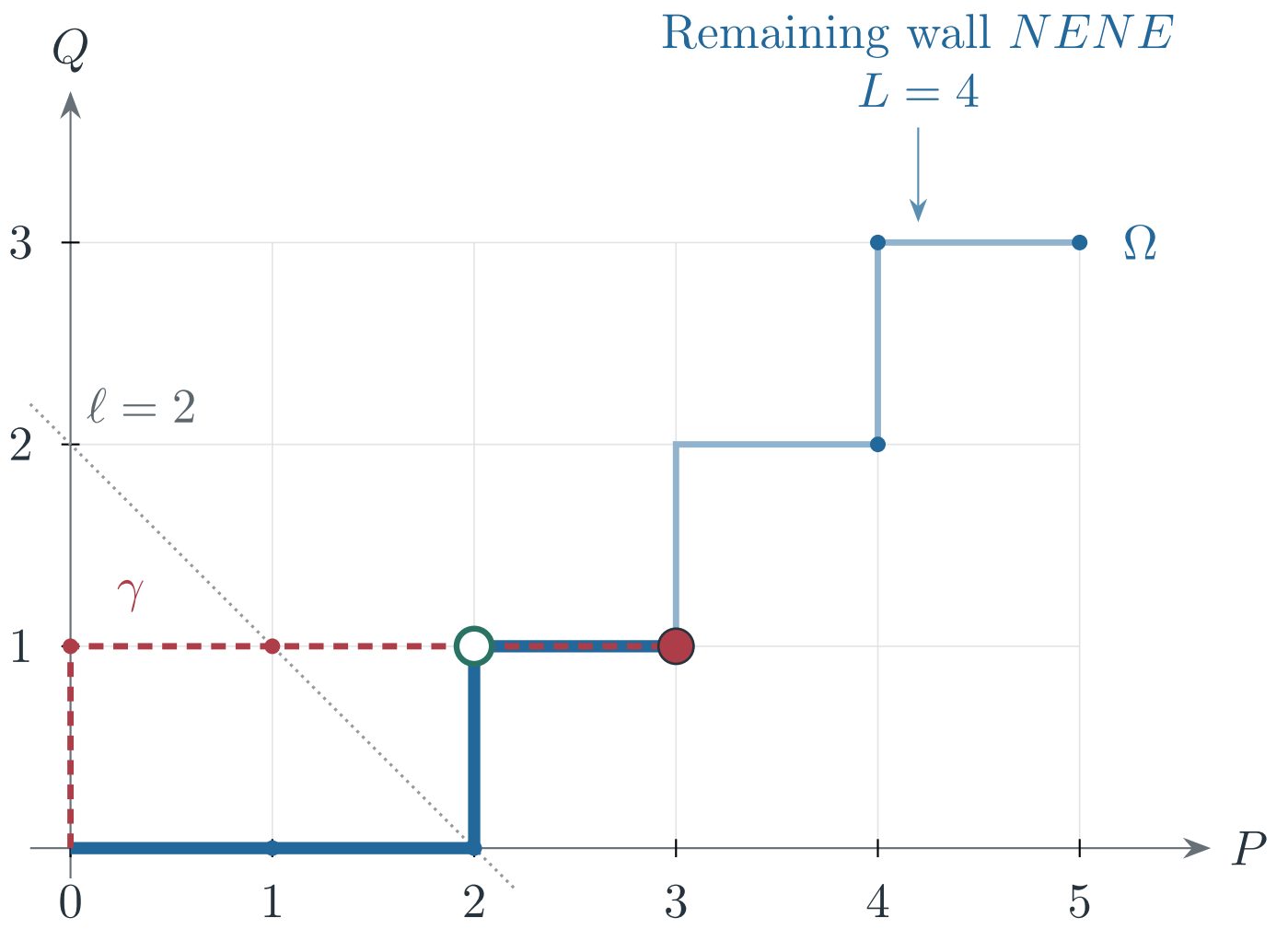
\vspace{1mm}
\centering
\begin{tikzpicture}[x=1.45cm,y=1.45cm,font=\small]
 \path[use as bounding box] (-.6,-.65) rectangle (6,4.25);
 \draw[step=1,draw=ink!13,line width=.3pt] (0,0) grid (5,3);
 \draw[axis] (-.2,0)--(5.65,0) node[right] {$P$};
 \draw[axis] (0,-.15)--(0,3.75) node[above] {$Q$};
 \foreach \i in {0,...,5} {\draw (\i,.045)--(\i,-.045) node[below=2pt] {\i};}
 \foreach \i in {1,...,3} {\draw (.045,\i)--(-.045,\i) node[left=2pt] {\i};}
 \draw[ink!50,densely dotted,line width=.6pt] (-.2,2.2)--(2.2,-.2);
 \node[anchor=south west,fill=white,inner sep=1pt,text=ink!75] at (.05,2.08) {$\ell=2$};
 \draw[draw=wallblue!50,line width=1.3pt] (3,1)--(3,2)--(4,2)--(4,3)--(5,3);
 \draw[wall,line width=2.5pt] (0,0)--(2,0)--(2,1)--(3,1);
 \draw[slack,line width=1.3pt] (0,0)--(0,1)--(3,1);
 \foreach \x/\y in {1/0,2/0,4/2,4/3,5/3} {\fill[wallblue] (\x,\y) circle[radius=1.6pt];}
 \fill[pathred] (0,1) circle[radius=1.6pt];
 \fill[pathred] (1,1) circle[radius=1.6pt];
 \draw[safeolive,line width=1.2pt,fill=white] (2,1) circle[radius=3.6pt];
 \filldraw[pathred,draw=ink,line width=.5pt] (3,1) circle[radius=3.6pt];
 \draw[draw=wallblue!75,-{Stealth[length=1.7mm]}] (4.2,3.57)--(4.2,3.1);
 \node[anchor=south,align=center,text=wallblue] at (4.2,3.58) {Remaining wall $NENE$\\$L=4$};
 \node[anchor=west,text=wallblue] at (5.12,3) {$\Omega$};
 \node[anchor=west,text=pathred,fill=white,inner sep=1pt] at (.2,1.25) {$\gamma$};
 \node[anchor=north,align=center] at (2,-.51) {$\Omega=EENENENE,\qquad\gamma=NEEE$};
\end{tikzpicture}

\vspace{3mm}
At $\ell=2$, the wall is at $(2,0)$ and the slack path at $(1,1)$;
the difference of their abscissae is $g=2-1=1$.
\end{minipage}
\par\vspace{5mm}

\begin{minipage}[t]{\linewidth}
\textbf{(b) The gap at successive diagonal levels}\par
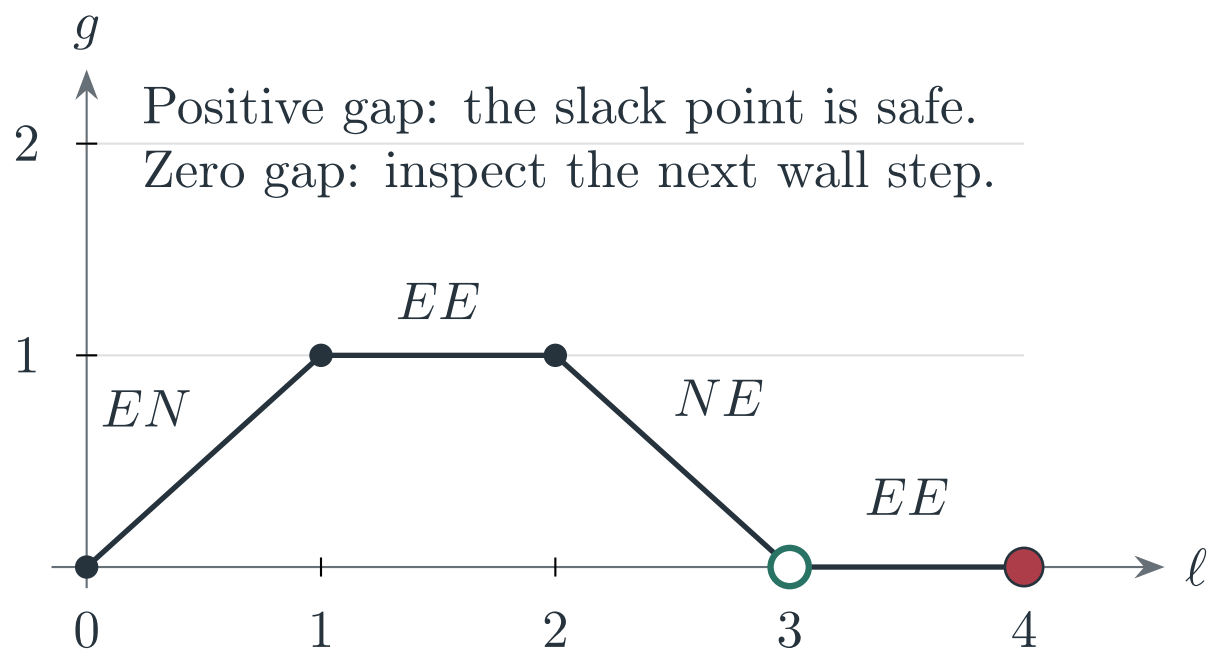
\vspace{1mm}
\centering
\begin{tikzpicture}[x=1.55cm,y=1.4cm,font=\small]
 \path[use as bounding box] (-.48,-1.1) rectangle (4.78,3.08);
 \draw[axis] (-.15,0)--(4.6,0) node[right] {$\ell$};
 \draw[axis] (0,-.1)--(0,2.35) node[above] {$g$};
 \draw[ink!15] (0,1)--(4,1);
 \draw[ink!15] (0,2)--(4,2);
 \foreach \i in {0,...,4} {\draw (\i,.045)--(\i,-.045) node[below=3pt] {\i};}
 \foreach \i in {1,2} {\draw (.045,\i)--(-.045,\i) node[left=3pt] {\i};}
 \draw[ink,line width=1pt] (0,0)--(1,1)--(2,1)--(3,0)--(4,0);
 \foreach \x/\y in {0/0,1/1,2/1} {\fill[ink] (\x,\y) circle[radius=2.2pt];}
 \draw[safeolive,line width=1.2pt,fill=white] (3,0) circle[radius=3.6pt];
 \filldraw[pathred,draw=ink,line width=.5pt] (4,0) circle[radius=3.6pt];
 \node[above left=1pt] at (.55,.55) {$EN$};
 \node[above=3pt] at (1.5,1) {$EE$};
 \node[above right=1pt] at (2.4,.6) {$NE$};
 \node[above=6pt] at (3.5,0) {$EE$};
 \node[anchor=north,align=center] at (2,-.62) {$g_\ell=P_\Omega(\ell)-P_\gamma(\ell)$\\$g_0,g_1,g_2,g_3,g_4=0,1,1,0,0$};
 \node[anchor=south west,align=left] at (.15,1.65) {Positive gap: the slack point is safe.\\Zero gap: inspect the next wall step.};
\end{tikzpicture}

\vspace{3mm}
\raggedright
$\circ$\quad At $(2,1)$, the next wall step is $E$: safe.\par
$\bullet$\quad At $(3,1)$, the next wall step is $N$: forbidden.
\end{minipage}
\par
\caption{\textbf{Safe contact and first forbidden contact.} A \emph{contact} is a common vertex of the two paths before the wall ends.
The slack path first returns to the wall at $(2,1)$, where the following east step
makes the contact safe.  The paired step $EE$ leads to $(3,1)$; this contact is
forbidden because the remaining wall starts north.  Thus first passage can occur
without the gap becoming negative.  The slack path stops at that contact, while
four wall steps remain.  This is a schematic example of the contact rule,
independent of the polyomino used in Figures~\ref{fig:core-shadow}--\ref{fig:fixed-boundary}. Here $P_\eta(\ell)$ denotes the abscissa of the path $\eta$ at diagonal level $\ell$. The TikZ code for this figure was prepared with assistance from ChatGPT (OpenAI) and reviewed by the author.}
\label{fig:contacts}
\end{minipage}
\end{figure}

A \emph{contact} is a common vertex of the two paths at equal diagonal
level, with positive remaining wall length.  It is safe or forbidden
according as its slack point is safe or forbidden.  A safe excursion
starts at a safe contact, ends at contact, and has positive gap at
every intermediate diagonal level.
By Corollary~\ref{pc:contact-rule}, a contact is safe exactly when
the next wall unit is east.
A safe excursion from contact either takes a paired $EE$ step or leaves by
$EN$ and first returns by $NE$.  On the two possible contact direction pairs $EE,NE$,
the resulting matrix is therefore
\[
 K_0=q\binom{\lambda}{a}(1,x\chi).
\]
At the first forbidden contact the remaining wall begins with a north unit.
If its length is $L$, horizontal overshoot and subsequent unrestricted tail
completion give $S_L=\mathcal T_1+\cdots+\mathcal T_L$.
Figure~\ref{fig:contacts} shows a safe contact followed by a first forbidden contact, with $L=4$ wall steps remaining. Its second panel gives the gap at successive diagonal levels.

For a series
$F(y)=\sum_{L\ge1}f_L(x)y^L$ whose coefficient of each power of $x$ is a
polynomial in $y$, \eqref{pe:tail-evaluations} gives
\begin{equation}\label{pe:suffix-functional}
 \mathcal S[F]:=\sum_{L\ge1}f_LS_L
 =\frac{F(\rho_+)-F(\rho_-)}{2\sqrt x}.
\end{equation}
Put $C_N(y)=C_N(y,y)$, so $y$ records total remaining length.  Then
\[
 C_N(y)=\frac{1-ry}{1-2\lambda y+(\lambda-x)y^2}.
\]
Since
\[
 yC_N(y)\bigm|_{y=(1\mp\sqrt x)^{-1}}
 =\frac{1-r\mp\sqrt x}{(1-\lambda)(1\mp2\sqrt x)},
\]
we find
\begin{equation}\label{pe:S-evaluation}
 \mathcal S[yC_N(y)]=\frac{1-2r}{(1-\lambda)(1-4x)}=:S.
\end{equation}
The completed north suffix has weights $bS,\lambda S$ after previous wall
directions $E,N$.  Starting from a returned contact of type $NE$, its total
weight $H_k$ is
\[
 H_k=(0,1)(I-K_0)^{-1}\binom{bS}{\lambda S}.
\]
Put $J=1-q\lambda-qa x\chi$.  The last row of the inverse is
$(qa,1-q\lambda)/J$, and
\[
 qab+\lambda-q\lambda^2=qx(1-\lambda),\qquad
 J=\frac{q(1-\lambda)}{1+t}.
\]
The first equality uses $ab=x-\lambda+\lambda^2$, and the second uses
$a\chi=t-(1+t)\lambda$.  Thus
\begin{equation}\label{pe:H}
 H_k=\frac{x(1+t)(1-2r)}{(1-4x)(1-\lambda)}.
\end{equation}
This includes every number of safe contacts and then the first forbidden north
suffix, with its complete unrestricted continuation.

\subsection{The initial boundary and the final sum}

The first canonical north slack step pairs with the compulsory first
horizontal wall unit and creates gap one.  The preceding wall direction
for the remainder counted by $C_N$ is north.  The required column,
indexed by the preceding slack direction $E,N$, is
\[
 \mathbf d_N=(d_{NE},d_{NN})^T=t\binom1{1+t}.
\]
By \eqref{pe:descent-powers}, every further initial horizontal wall unit
is represented by multiplication by
\[
 K_t=V\begin{pmatrix}1&0\\0&t\end{pmatrix}.
\]
In the flat case the forced north wall unit is represented by multiplication by
\[
 N_t=bV\begin{pmatrix}t^{-1}&0\\0&1\end{pmatrix}.
\]
The factor $t^{-1}$ removes one already counted descent when the slack step
is east; even when it is the first return to zero, $t^{-1}d_{NE}=1$ gives
the correct immediate contact.  Because the initial horizontal part has at
least one unit, this step never starts at negative gap.

For $P(w)=\sum_{m\ge1}p_m(x)w^m$, write $\widetilde P(w)=P(w)/w$
to remove the first east step already paired with the canonical step.
The matrix substitution used below means
\[
 \widetilde P(K)=\sum_{m\ge1}p_m(x)K^{m-1}.
\]
These sums are coefficientwise finite: at fixed semi-perimeter the length
of the distinguished initial segment is bounded.  Let $L_k^N,L_k^F$
count the excluded ascending polyominoes with core degree $k$ in the
nonflat class $N_k^{\mathrm{lo}}$ and the flat class
$F_{k+1}^{\mathrm{av}}$, respectively, including every column after the
first forbidden one.  Equations \eqref{pe:ordered-walls} and
\eqref{pe:north-suffix} give
\[
 L_k^N=xZ_0H_k(0,1)\widetilde P_N(K_t)\mathbf d_N,
 \qquad
 L_k^F=xH_k(0,1)\widetilde P_F(K_t)N_t\mathbf d_N.
\]
These matrix expressions reduce to scalar evaluations.  Indeed
\[
 K_t\mathbf d_N=(1+t)\mathbf d_N,
 \qquad N_t\mathbf d_N=\frac{b(1+t)}t\mathbf d_N,
\]
by $x=t/(1+t)^2$.  Therefore
\begin{equation}\label{pe:scalar-boundary}
 L_k^N=xZ_0H_k\,tP_N(1+t),\qquad
 L_k^F=xH_k\,b(1+t)P_F(1+t).
\end{equation}
The subfamily $W_N(w;u,0)$ has no forbidden north suffix and so contributes zero,
as required for the exact non-flat class.

Put $z=t^k$ and $h=1+t$.  The quantities needed to evaluate
\eqref{pe:scalar-boundary} are
\begin{align*}
 r&=\frac{t-z}{h(1-z)},&
 s&=\frac{t(1-z)}{h(1-tz)},&
 B&=\frac{t(1-tz)}{h(1-t^2z)},\\
 \varrho&=\frac{t^2-z}{h(t-z)},&
 \lambda&=\frac{t(1-z^2)}{h(1-tz^2)},&
 Z_0&=\frac{h(1-z)(1-tz)}{(1-t)(1-tz^2)},\\
 a&=\frac{z(1-t)(1-tz)}{h(1-z)(1-tz^2)},&
 b&=\frac{tz(1-t)(1-z)}{h(1-tz)(1-tz^2)},&
 \chi&=\frac{tz h(1-z)}{1-tz},\\
 S_N&=\frac{tz(1-t)^2}{h^2(t-z)(1-tz)},&
 S_F&=\frac{tz(1-t)^2}{h^2(1-z)(1-t^2z)}.&&
\end{align*}
Substitution in \eqref{pe:PN}--\eqref{pe:PF}, using
\eqref{pe:tree-marking}, gives the two scalar identities
\begin{equation}\label{pe:boundary-evaluations}
 P_N(h)=\frac{t^2z}{h^4(1-tz)},\qquad
 P_F(h)=\frac{t^3z(1-tz)^2}{(1-t)h^3(1-t^3z^2)}.
\end{equation}
Since $xh=t/h$ is a fixed point of $f$, one has $T_j(h)=t/h$ for
every $j\ge1$.  The remaining factors reduce to
\[
 1-\varrho h=\frac{t(1-t)}{t-z},\quad
 1-T_{k-1}(h)-T_k(h)=\frac{1-t}{h},\quad
 1-rh=\frac{1-t}{1-z},
\]
\[
 1-T_k(h)-s=\frac{1-t}{h(1-tz)},\qquad
 1-B-s=\frac{(1-t)(1-t^3z^2)}{h(1-tz)(1-t^2z)}.
\]
These identities together with
\eqref{pe:H} give
\begin{align}
 L_k^N&=\frac{t^5z(1+z)}{(1-t)^2(1+t)^4(1-t^2z^2)},
 \label{pe:LN}\\
 L_k^F&=\frac{t^6z^2(1+z)}
 {(1-t)(1+t)^4(1+tz)(1-t^3z^2)}.
 \label{pe:LF}
\end{align}
Thus, writing $Z=t^{k+1}=tz$ and combining denominators,
\begin{equation}\label{pe:later-total}
 L_k^N+L_k^F
 =\frac{t^3Z(t+Z)(1+Z-Z^2-tZ)}
 {(1-t)^2(1+t)^4(1-Z^2)(1-tZ^2)}.
\end{equation}
The first-passage decomposition \eqref{thm:first-passage-factorisation}
and \eqref{pe:bad-core} now yield
\begin{equation}\label{pe:excluded}
 \boxed{
 A(x)-A_k(x)=
 \frac{t^{2k+4}(1+3t^2)+2t^{k+4}(1+t)}
 {(1-t)^3(1+t)^5(1-t^{2k+2})}.}
\end{equation}
Indeed adding \eqref{pe:later-total} to $R(t,Z)$ cancels $1-tZ^2$ and
leaves numerator $t^2Z^2(1+3t^2)+2t^3Z(1+t)$ over
$(1-t)^3(1+t)^5(1-Z^2)$, which is \eqref{pe:excluded}.

\section{Proof of the generating function}\label{sec:conclusion}
The unrestricted convex series is the classical formula of
Delest and Viennot~\cite[Theorem~11.1]{DV}.  With $x$ marking semi-perimeter,
it is
\begin{equation}\label{eq:classical-convex}
 C(x)=\frac{x^2(1-6x+11x^2-4x^3)}{(1-4x)^2}
                  -\frac{4x^4}{(1-4x)^{3/2}}.
\end{equation}

The preceding decomposition divides $A-A_k$ into cores whose degree is
already greater than $k$ and the first degree increase in a tail attached to
a core of degree $k$.  Their sum is \eqref{pe:excluded}, with $x=t/(1+t)^2$.
The directional decomposition and \eqref{eq:classical-convex} therefore give
$Q_k=C-2(A-A_k)$ explicitly.

To put this expression in the form of Theorem~\ref{thm:main}, solve the two
polynomial recurrences after the same substitution:
\[
 d_k=\frac{1-t^{2k+2}}{(1-t^2)(1+t)^{2k}},\qquad
 r_k=\frac{1+t^k}{(1+t)^k},\qquad
 \sqrt{1-4x}=\frac{1-t}{1+t}.
\]
These identities follow at once by inserting the two roots
$1/(1+t)^2,t^2/(1+t)^2$ of the first characteristic equation and
$1/(1+t),t/(1+t)$ of the second, including the stated initial values.
Substituting them and the definition of $p_k$ in \eqref{eq:main}, subtracting
from \eqref{eq:classical-convex}, and dividing by two gives exactly
\eqref{pe:excluded}.  This proves Theorem~\ref{thm:main}.

\section*{Funding}
This research did not receive any specific grant from funding agencies in
the public, commercial, or not-for-profit sectors.

\section*{Declaration of competing interests}
The author declares that he has no known competing financial interests or
personal relationships that could have appeared to influence the work
reported in this paper.

\section*{Data availability}
This article presents theoretical results. The definitions, constructions,
and proofs are included in the article; no external dataset is required
to establish the results.

\section*{Declaration of generative AI and AI-assisted technologies
in the manuscript preparation process}
During the preparation of this work, the author used ChatGPT (OpenAI)
for mathematical exploration and preparation of the TikZ figures,
Claude (Anthropic) for critical feedback, and Codex
(OpenAI) for \LaTeX{} formatting and preparation of the submission files.
After using these tools, the author reviewed and edited the content as
needed and takes full responsibility for the content of the published
article.

\end{document}